%% file: KaKuBr2017rspa.tex
\documentclass[openacc]{rsproca_new}%%%%where rsproca is the template name

\input{Header_Pkgs_Defs}

\newcommand{\RefOne}[1]{{\color{magenta} \textbf{#1}}}
\newcommand{\RefTwo}[1]{{\color{blue} \textbf{#1}}}
\begin{document}

%%%% Article title to be placed here
\title{Sparse identification of nonlinear dynamics for model predictive control in the low-data limit}

\author{%%%% Author details
E. Kaiser$^{1}$, J.~N. Kutz$^{2}$ and S.~L. Brunton$^{1}$}

%%%%%%%%% Insert author address here
\address{$^{1}$Department 
of Mechanical Engineering, University of Washington, Seattle, WA, 98195\\
$^{2}$Department
of Applied Mathematics, University of Washington, Seattle,
WA, 98195}

%%%% Subject entries to be placed here %%%%
\subject{xxxxx, xxxxx, xxxx}

%%%% Keyword entries to be placed here %%%%

%%%% Insert corresponding author and its email address}
\corres{Eurika Kaiser\\
\email{eurika@uw.edu}}

%%%% Abstract text to be placed here %%%%%%%%%%%%
\begin{abstract}
\input{Abstract}

\end{abstract}
%%%%%%%%%%%%%%%%%%%%%%%%%%%

%%%%%%%%%% Insert the texts which can accomdate on firstpage in the tag "fmtext" %%%%%

\begin{fmtext}

\end{fmtext}

%%%%%%%%%%%%%%% End of first page %%%%%%%%%%%%%%%%%%%%%

\maketitle

%%%%%%%%%%%%
%%% INTRODUCTION
%%%%%%%%%%%%
\input{Sec1}

%%%%%%%%%%%%
%%% BACKGROUND
%%%%%%%%%%%%
\input{Sec3}

\input{Sec41}
\input{Sec42}
\input{Sec43}

\input{Sec44}

\input{Sec5}

%[!h]

\enlargethispage{20pt}

\ethics{No ethical considerations apply.}

\dataccess{The code used in this work is made available at:\\ \url{https://github.com/eurika-kaiser/SINDY-MPC}. The data can be generated using the code.}

\competing{We declare we have no competing interests.}

\aucontribute{All authors conceived of the work, designed the study and drafted the manuscript. EK carried out the computations.}

\ack{The authors gratefully acknowledge many valuable discussions with Josh Proctor. }
%We would also like to acknowledge ....
%\disclaimer{Insert disclaimer text here.}

\funding{
	EK gratefully acknowledges support by the Washington Research Foundation, the Gordon and Betty Moore Foundation (Award \#2013-10-29), the Alfred P. Sloan Foundation (Award \#3835), and the University of Washington eScience Institute.
	SLB and JNK acknowledge support from the Defense Advanced Research Projects Agency (DARPA contract HR011-16-C-0016 and PA-18-01-FP-125). SLB acknowledges support from the Army Research Office (W911NF-17-1-0306 and W911NF-17-1-0422). JNK acknowledges support from the Air Force Office of Scientific Research (FA9550-17-1-0329).}

%EK gratefully acknowledges support by the ``Washington Research Foundation Fund for Innovation in Data-Intensive Discovery" and a Data Science Environments project award from the Gordon and Betty Moore Foundation (Award \#2013-10-29) and the Alfred P. Sloan Foundation (Award \#3835) to the University of Washington eScience Institute.

%%%%%%%%%% Insert bibliography here %%%%%%%%%%%%%%

%\vskip2pc

\bibliographystyle{RS}
\bibliography{references}

%\noindent {\bf Please follow the coding for references as shown below.}

%\begin{thebibliography}{9}
%
%\bibitem{1} Allwood JM, Cullen JM. 2011 \textit{Sustainable materials:  with both eyes open}.
%Cambridge, UK: UIT Cambridge. See \href{http://www.withbotheyesopen.com}{http://www.withbotheyesopen.com}.
%
%\bibitem{2}  MacKay DJC. 2008  \textit{Sustainable energy:  without the hot air}.
% Cambridge, UK: UIT Cambridge. See \href{http://www.withouthotair.com}{http://www.withouthotair.com}.
%
%\bibitem{3} Gallman PG. 2011  \textit{Green alternatives and national energy strategy: the facts
% behind the headlines}.  Baltimore,\ MD: Johns Hopkins University Press.
%
%\bibitem{4} MacKay DJC. 2013.  Solar energy in the context of energy use, energy transportation, and
% energy storage. \textit{Proc. R. Soc. A} \textbf{371}.
%
%\end{thebibliography}

%\noindent If maintaining .bib file for references, then please use "RS.bst" to generate the references.

%\noindent Example:

%\verb+\bibliographystyle{RS}+ %%%% .BST file

%\verb+\bibliography{references}+ %%%%% .Bib file

\end{document}

%% file: Header_Pkgs_Defs
\usepackage{amsmath, amssymb} %amsthm
\usepackage{graphicx}
\usepackage{epstopdf}
\usepackage{overpic}
\usepackage{cancel}
\usepackage{rotating}
\usepackage{url}
\usepackage[usenames,dvipsnames,svgnames,table]{xcolor}
\usepackage{rotating}
\usepackage{multirow}
\usepackage{multicol}
\usepackage{multirow}
\usepackage{adjustbox}
\usepackage{setspace}
\usepackage{algorithm,algcompatible}
\usepackage{algpseudocode}
\usepackage{tikz}
\usepackage{booktabs}
\renewcommand{\Comment}[2][.4\linewidth]{%
	\leavevmode\hfill\makebox[#1][l]{$\triangleright$~#2}}
\usepackage{mathtools}
\usepackage{psfrag}
\usepackage{tikz}
\usepackage[T1]{fontenc}
\usepackage{array}
\usepackage{makecell}
\newcolumntype{x}[1]{>{\centering\arraybackslash}p{#1}}

\DeclareMathOperator*{\argmin}{arg\rm{}min}
\algnewcommand\algorithmicinput{\textbf{Input:}}
\algnewcommand\INPUT{\item[\algorithmicinput]}
\algnewcommand\algorithmicoutput{\textbf{Output:}}
\algnewcommand\OUTPUT{\item[\algorithmicoutput]}
\algnewcommand\algorithmicoptional{\textbf{Optional:}}
\algnewcommand\OPTIONAL{\item[\algorithmicoptional]}

\newcommand{\bF}{\mathbf{F}}

\newcommand{\bI}{\mathbf{I}}
\newcommand{\bK}{\mathbf{K}}

\newcommand{\bQ}{\mathbf{Q}}

\newcommand{\bR}{\mathbf{R}}

\newcommand{\bu}{\mathbf{u}}

\newcommand{\bX}{\mathbf{X}}
\newcommand{\bx}{\mathbf{x}}

\newcommand{\by}{\mathbf{y}}

\definecolor{blue}{rgb}{0,0,1}
\definecolor{darkgreen}{rgb}{0,0.5,0}
\definecolor{red}{rgb}{1,0,0}
\definecolor{teal}{rgb}{0,0.5,0.7}

\newcommand{\boundellipse}[3]% center, xdim, ydim
{(#1) ellipse (#2 and #3)
}

%% file: Abstract.tex
%The data-driven discovery of dynamics via machine learning is currently pushing the frontiers of modeling and control efforts, and it provides a tremendous opportunity to extend model predictive control to systems that do not admit first-principles models.  
Data-driven discovery of dynamics via machine learning is pushing the frontiers of modeling and control efforts, providing a tremendous opportunity to extend the reach of model predictive control.  
However, many leading methods in machine learning, such as neural networks, require large volumes of training data, may not be  interpretable, do not easily include known constraints and symmetries, and may not generalize beyond the attractor where models are trained.  
These factors limit their use for the online identification of a model in the low-data limit, for example following an abrupt change to the system dynamics.
%These factors limit the use of these techniques for the online identification of a model in the low-data limit, for example following an abrupt change to the system dynamics.  
%
In this work, we extend the recent sparse identification of nonlinear dynamics (SINDY) modeling procedure to include the effects of actuation and demonstrate the ability of these models to enhance the performance of model predictive control (MPC), based on limited, noisy data.  
%and demonstrate the ability of these models to provide enhanced performance, based on limited data, when combined with model predictive control (MPC).  
%
SINDY models are parsimonious, identifying the fewest terms in the model needed to explain the data, making them interpretable and generalizable. %, and reducing the burden of training data. 
We show that the resulting SINDY-MPC framework has higher performance, requires significantly less data, and is more computationally efficient and robust to noise than neural network models, making it viable for online training and execution in response to rapid system changes. %changes to the system.  
SINDY-MPC also shows improved performance over linear data-driven models, although linear models may provide a stopgap until enough data is available for SINDY.  
SINDY-MPC is demonstrated on a variety of dynamical systems with different challenges, including the chaotic Lorenz system, a simple model for flight control of an F8 aircraft, and an HIV model incorporating drug treatment. %\vspace{-.2in}

%% file: Sec1.tex
\section{Introduction}
%%%%
%%%% PARAGRAPH 1
%%%%
Data-fueled modeling and control of complex systems is currently undergoing a revolution, driven by the confluence of big data, advanced algorithms in machine learning, and modern computational hardware.  
Model-based control strategies, such as model predictive control, are ubiquitous, relying on accurate and efficient models that capture the relevant dynamics for a given objective.  
Increasingly, first principles models are giving way to data-driven approaches, for example in turbulence, epidemiology, neuroscience, and finance~\cite{Kutz2016book}.
Although these methods offer tremendous promise, there has been slow progress in distilling physical models of dynamic processes from data. 
%
%% ALTERNATE ORDER
% Model-based control strategies, such as model predictive control, are ubiquitous, relying on accurate and efficient models that capture the relevant dynamics for a given objective.  
%%
% Increasingly, first principles models are giving way to data-driven approaches, for example in  turbulence, epidemiology, neuroscience, and finance~\cite{Kutz2016book}. 
%%
% Indeed, the data-driven modeling and control of complex systems is currently undergoing a revolution, driven by the confluence of big data, advanced algorithms in machine learning, and modern computational hardware.  
%%
% Although these methods offer tremendous promise, there has been slow progress in distilling physical models of dynamic process from data. 
%%
Despite their undeniable success, many modern techniques in machine learning (e.g., neural networks) rely on access to massive data sets, have limited ability to generalize beyond the attractor where data is collected, and do not readily incorporate known physical constraints.  
The current challenges associated with data-driven discovery limit its use for real-time control of strongly nonlinear, high-dimensional, multi-scale systems, and prevent online recovery in response to abrupt changes in the dynamics. 
Fortunately, a new paradigm of sparse and parsimonious modeling is enabling interpretable models in the low-data limit.
In this work, we extend the recent sparse identification of nonlinear dynamics (SINDy) framework~\cite{Brunton2016pnas} to identify models with actuation, and combine it with model predictive control (MPC) for effective and interpretable data-driven, model-based control. 
We apply the proposed SINDY-MPC architecture to control several nonlinear systems and demonstrate improved control performance in the low-data limit, compared with other leading data-driven methods, including linear response models and neural networks.  

%\begin{center}
%	\it Paragraph 1: 
%\end{center}
%\begin{itemize}\setlength{\itemsep}{0pt}
%	\item Similar to Para 1 in SINDYc paper
%	\item Recent advances in data-driven discovery of dynamical systems are transformative
%	\item New opportunities for model-based control strategies  (LQR, SDRE, MPC (in particular), predictive control?): provide physical insights, probably better performance as a better representation of the true underlying dynamics, low-dimensional
%\end{itemize}

% Point of paragraph 1: 
% Data-driven modeling is exploding
% Model-based methods increasingly relying on data as first principles are failing us. 
% Data-intensive methods, such as neural networks? require significant amounts of data and offline training.  Move towards low-data limit and online characterization.  Physical interpretability, addition of constraints
% A new paradigm of sparse/parsimonious models are poised to capitalize on low-data limit. 
% In this paper, we combine a recently developed method with MPC for real-time, online characterization of systems.  We compare with linear models (ERA/DMDc/impulse) and NARX models.  Kick ass.  

%%%%
%%%% PARAGRAPH 2
%%%%
Model-based control techniques, such as MPC~\cite{allgower1999springer,camacho2013model} and optimal control~\cite{sp:book,dp:book}, are cornerstones of advanced process control, and are well-positioned to take advantage of the data-driven revolution.  
Model predictive control is particularly ubiquitous in industrial applications, as it enables the control of strongly nonlinear systems with constraints, which are difficult to handle using traditional linear control approaches~\cite{garcia1989model,morari1999model,lee2011springer,mayne2014automatica,Utku2017jgcd}.   %qin1997proc,rawlings2000tutorial,qin2003survey,garriga2010model
MPC benefits from simple and intuitive tuning and the ability to control a range of simple and complex phenomena, including systems with time delays, non-minimum phase dynamics, and instability.  
In addition, it is straightforward to incorporate known constraints and multiple operating conditions, it exhibits an intrinsic compensation for dead time, and it provides the flexibility to formulate and tailor a control objective.
The major drawback of model-based control, such as MPC, lies in the development of a suitable model via existing system identification or model reduction techniques~\cite{Brunton2015amr}, which may require expensive and time-consuming data collection and computations. 

Nearly all industrial applications of MPC rely on empirical models, and increasing plant complexity and tighter performance specifications require models with higher accuracy.  
There are many techniques to obtain data-driven models, including state-space models from the eigensystem realization algorithm (ERA)~\cite{ERA:1985} and other subspace identification methods, Volterra series~\cite{Brockett1976automatica,Boyd1984imajmci,maner1994nonlinear}, autoregressive models~\cite{Akaike1969annals} (e.g., ARX, ARMA, NARX, and NARMAX~\cite{Billings2013book} models), and neural network models~\cite{lippmann1987introduction,draeger1995model,wang2016combined,aggelogiannaki2008nonlinear}, to name only a few. 
These procedures all tend to yield black-box models, with limited interpretability, physical insights, and ability to generalize.  
More recently, linear representations of nonlinear systems using extended dynamic mode decomposition~\cite{Williams2015jnls} have been successfully paired with MPC~\cite{korda2016_a,Peitz2017arxiv}. 
Nonlinear models based on machine learning, such as neural networks, are increasingly used due to advances in computing power, and recently deep reinforcement learning has been combined with MPC~\cite{Peng.2009,Zhang2016icra}, yielding impressive results in the large-data limit.  
However, large volumes of data are often a luxury, and many systems must be identified and controlled with limited data, for example in response to abrupt changes.  
Current efforts are focused on \emph{rapid} learning based on minimal data.   

\begin{figure}[H]
	\centering
	\begin{overpic}[width=\textwidth]{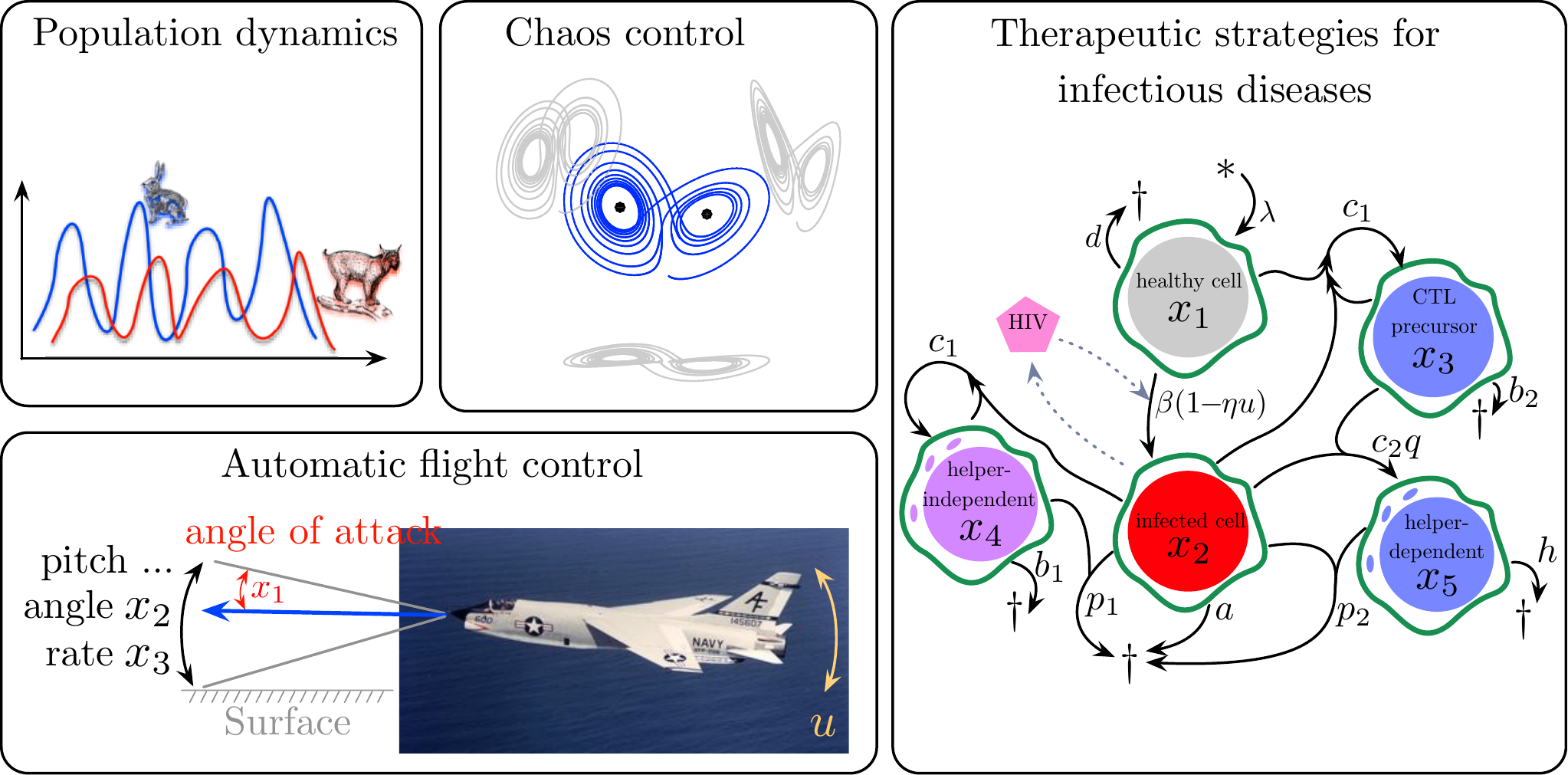}
		\put(10,44){(Sec.~\ref{Sec:Lotka-Volterra})}
		\put(37,44){(Sec.~\ref{Sec:Lorenz})}
		\put(43,19){(Sec.~\ref{Sec:F8-Aircraft})}
		\put(74,40){(Sec.~\ref{Sec:HIV})}			
	\end{overpic}	
	\vspace{-.2in}
	\caption{Applications of SINDY-MPC investigated in this work.}
	\vspace{-.2in}
	\label{Fig:OVERVIEW_APPS}
\end{figure}
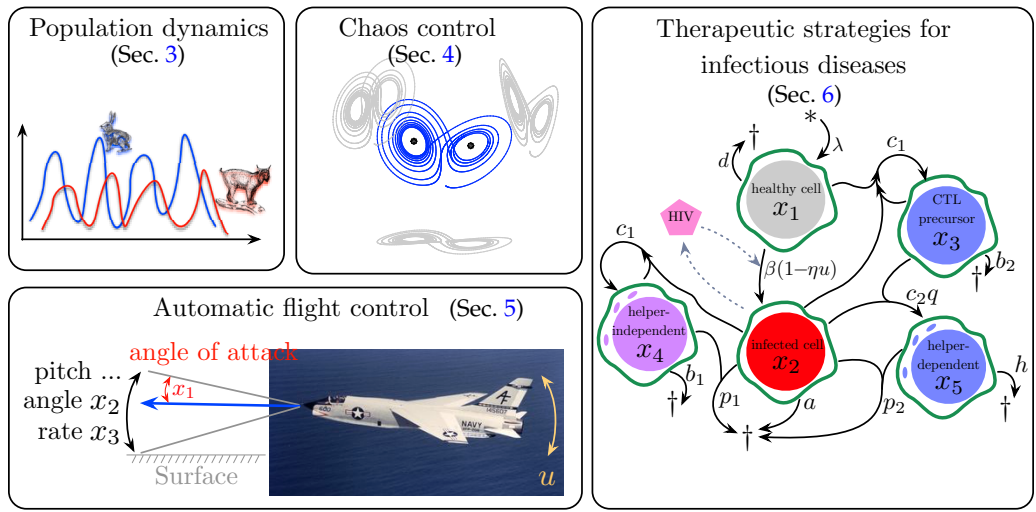
When abrupt changes occur in the system, an effective controller must rapidly characterize and compensate for the new dynamics, leaving little time for discovery based on limited data.  
A second challenge is the ability of models to generalize beyond the training data, which is related to the ability to incorporate new information and quickly modify the model.
Machine learning algorithms often suffer from overfitting and a lack of interpretability, although the application of these algorithms to physical systems offers a unique opportunity to incorporate known symmetries and constraints.  
These challenges point to the need for \emph{parsimonious} and interpretable models~\cite{Bongard2007pnas,Schmidt2009science,Brunton2016pnas} that may be characterized from limited data and in response to \RefTwo{abrupt changes~\cite{Quade2018arxiv}}.  
Whereas traditional methods require unrealistic amounts of training data, the recently proposed SINDY framework~\cite{Brunton2016pnas} relies on sparsity-promoting optimization to identify parsimonious models from limited data, resulting in interpretable models that avoid overfitting.  
It has also been shown recently~\cite{Loiseau2016jfm} that it is possible to enforce known physics (e.g., constraints, conservation laws, and symmetries) in the SINDY algorithm, improving stability and performance of models.  

%
%\vspace{0.5cm}
%\begin{center}
%	\it Paragraph 3: Data-driven models for control
%\end{center}
%\begin{itemize}\setlength{\itemsep}{0pt}
%	\item In the context of model complexity, linear/nonlinear, amount of training data, hyperparametrization, training time, execution time, robustness w.r.t. noise, parameter sensitivities, etc.
%	\item System identification (Para 2 in SINDYc paper + NN etc.) 
%	\item New avenues: DMD, Koopman, symbolic regression, SINDy, etc.
%	\item Sparsity discussion
%\end{itemize}

In this work, we combine SINDY with MPC for enhanced data-driven control of nonlinear systems in the low-data limit.  
First, we extend the SINDY architecture to identify interpretable models that include nonlinear dynamics and the effect of actuation.  
Next, we show the enhanced performance of SINDY-MPC compared with linear data-driven models and with neural network models.  
The linear models are identified using dynamic mode decomposition with control (DMDc)~\cite{Proctor2016siads,Kutz2016book}, which is closely related to SINDY and traditional state-space modeling techniques such as ERA. 
SINDY-MPC is shown to have better prediction accuracy and control performance than neural network models, especially for small and moderate amounts of noisy data.  
In addition, SINDY models are less expensive to train and execute than neural network models, enabling real-time applications.  
SINDY-MPC also outperforms linear models for moderate amounts of data, although DMDc provides a working model in the extremely low-data limit for simple problems.
Thus, in response to abrupt changes, a linear DMDc model may be used until a more accurate SINDY model is trained.  
We demonstrate the SINDY-MPC architecture on several systems of increasing complexity as illustrated in Fig.~\ref{Fig:OVERVIEW_APPS}.

%
%
%
%\vspace{0.5cm}
%\begin{center}
%	\it Paragraph 4: SINDYc for model-based control / MPC
%\end{center}
%\begin{itemize}\setlength{\itemsep}{0pt}
%	\item Extension of SINDY
%	\item Pro: Interpretability (in contrast to black-box models such as NNs), easier/faster training compared with neural networks, 
%	\item Put in discussion or mention briefly here: usable with high-dimensional data via SVD, time-delay coordinates (combined with sparse sampling for time delay selection), input--output data 
%\end{itemize}

%% file: Sec3.tex
%rojas2014ieee,rojas2013eec

\section{SINDY-MPC framework}
\begin{figure}[tb]
	%	\vspace{-0.2in}
	\centering
	\begin{overpic}[width=\textwidth]{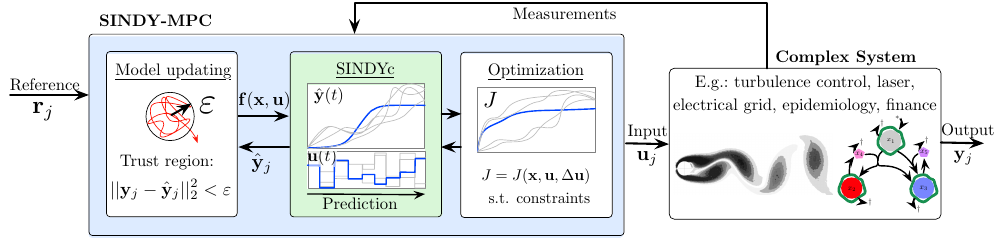}
	\end{overpic}	
	\vspace{-.2in}
	\caption{Schematic of the proposed SINDY-MPC framework, using sparse nonlinear models for predictive control.}\label{Fig:SINDYMPC}
		\vspace{-0.2in}
\end{figure}
The SINDY-MPC architecture combines the systematic data-driven discovery of dynamics with advanced model-based control to facilitate rapid model learning and control of strongly nonlinear systems.  
The overarching SINDY-MPC framework is illustrated in Fig.~\ref{Fig:SINDYMPC}. 
In the following sections, we will describe the sparse identification of nonlinear dynamics with control and model predictive control algorithms.  
We consider the nonlinear dynamical system
\begin{align}\label{Eq:NonlinearDynamicsWithControl}
\frac{d}{dt}\bx = {\bf f}(\bx,\bu),\quad \bx(0) = \bx_0
\end{align}
with state $\bx\in\mathbb{R}^n$, control input $\bu\in\mathbb{R}^q$, and smooth dynamics ${\bf f}(\bx,\bu): \mathbb{R}^n\times \mathbb{R}^q \rightarrow \mathbb{R}^n$.

%\section{Sparse identification of nonlinear dynamics with control (SINDYc) }

\subsection{Sparse identification of nonlinear dynamics with control}
Advanced machine learning algorithms provide new opportunities for nonlinear system identification. 
In particular, sparsity-promoting methods are playing an increasingly important role by recognizing the importance of parsimony in models~\cite{li2005ieee,chen2014ieee,Brunton2016pnas}, i.e. the tradeoff between model complexity and data fit. 
Recent work based on compressed sensing has been used to handle noise and outliers~\cite{xu2014automatica} \RefTwo{for linear system identification} and large libraries of candidate functions~\cite{pan2012cdc}. 
Sparse regularization, which has been demonstrated for parameter and structure identification~\cite{chen2014ieee,calafiore2015automatica,pan2016ieee,Brunton2016pnas}, is a particularly promising direction as this can promote robustness and generalizability in models. 
%In this context, compressed sensing based approaches have been shown to robustly identify linear systems in the presence of noise and outliers~\cite{xu2014automatica} and nonlinear systems using large libraries of candidate functions in a Bayesian setting~\cite{pan2012cdc}. 
%Sparse regularization, in particular, is a promising direction as this can promote robustness and generalizability in models. Recent efforts demonstrate its applicability for finite-impulse response models using kernel methods~\cite{chen2014ieee}, large-scale parameter identification~\cite{calafiore2015automatica}, and structure detection using Bayesian methods~\cite{pan2016ieee}. 
We refer the reader to an extensive review on nonlinear system identification methods~\cite{nelles2013springer} and a recent review in the context of machine learning~\cite{pillonetto2014automatica}.

\begin{figure}[tb]
	\centering
	\begin{overpic}[width=\textwidth]{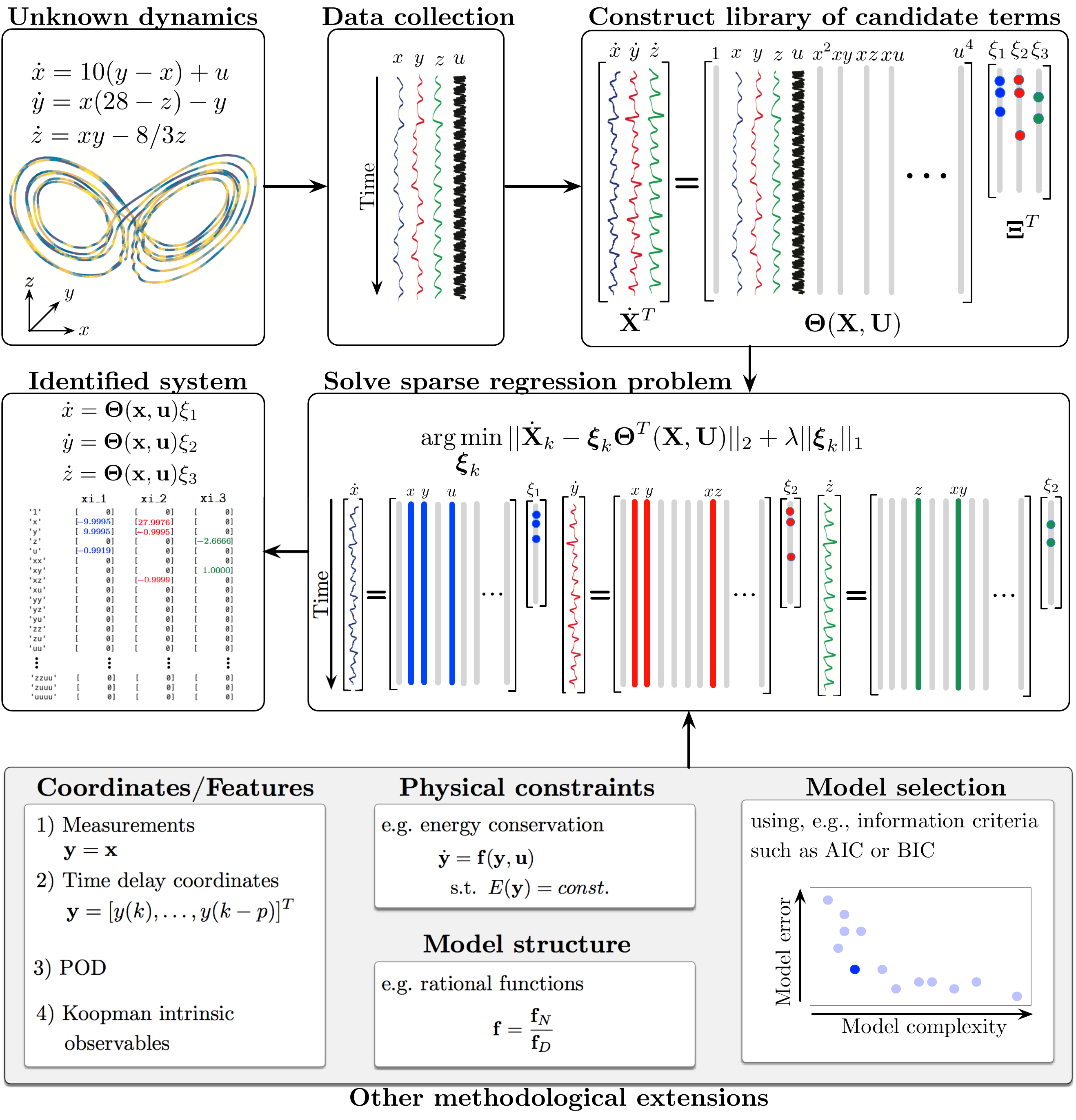}
		\put(25,21){\scriptsize \cite{Brunton2017natcomm}}
		\put(10,12){\scriptsize \cite{Kutz2016book}}
		\put(16,6){\scriptsize \cite{Kaiser2017arxiv}}
		\put(54.5,25.5){\scriptsize \cite{Loiseau2016jfm}}
		\put(53.5,11.5){\scriptsize \cite{Mangan2016ieee}}
		\put(84.5,24){\scriptsize \cite{Mangan2017prsa}}
	\end{overpic}
	\vspace{-.25in}
	\caption{Schematic of the SINDYc algorithm and extensions.  Active terms in a library of candidate nonlinearities are selected via sparse regression. 
	Illustration of the modular nature of the SINDY with control framework (bottom row) and its ability to handle high-dimensional systems, limited measurements, known physical constraints, and model selection.}
	\label{Fig:SINDYc:Schematic}	
	\vspace{-.25in}
\end{figure}
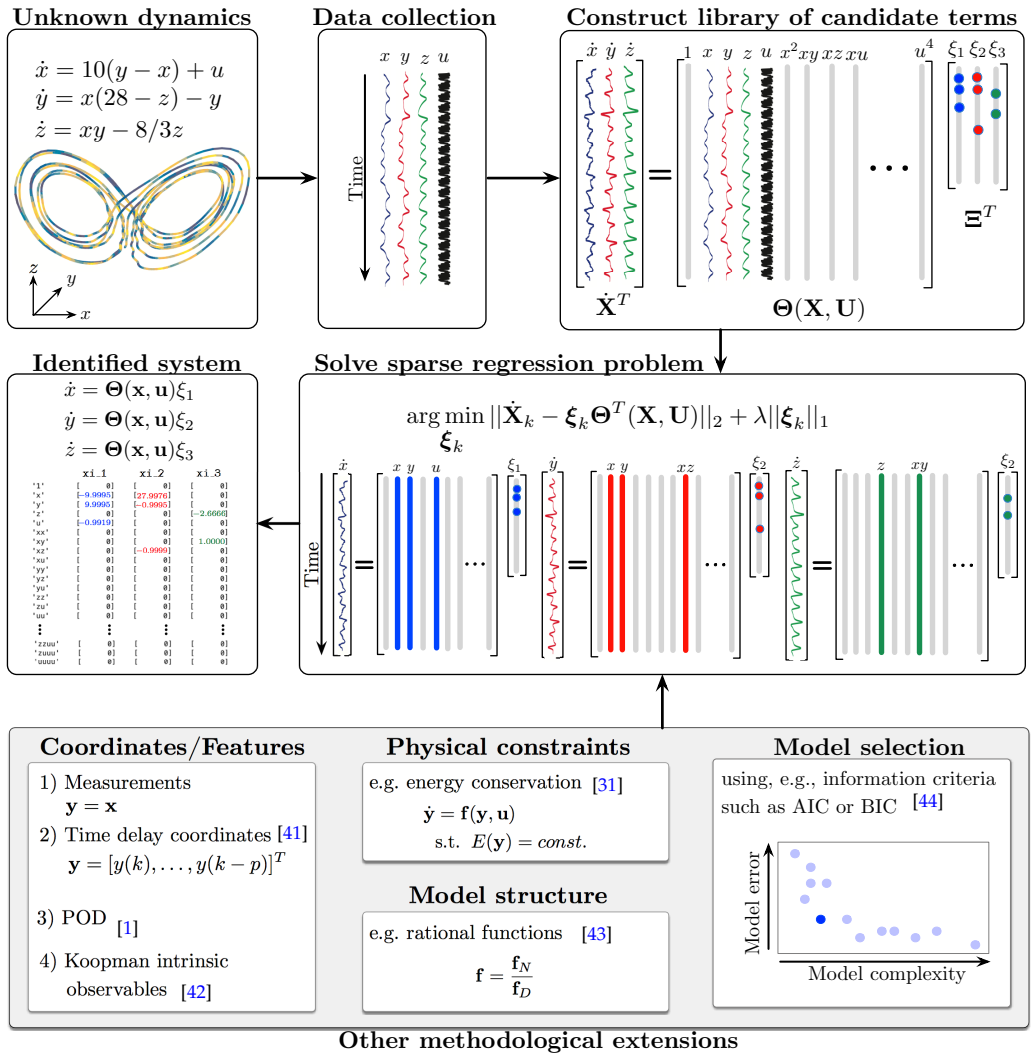

Here, we generalize the sparse identification of nonlinear dynamics (SINDY) method~\cite{Brunton2016pnas} to include inputs and control %, building on a previous conference paper~\cite{Brunton2016nolcos} and 
as illustrated in Fig.~\ref{Fig:SINDYc:Schematic}. 
SINDY identifies nonlinear dynamical systems from measurement data, relying on the fact that many systems have relatively few terms in the governing equations.  
Thus, sparsity-promoting techniques may be used to find models that automatically balance sparsity in the number of model terms with accuracy, resulting in parsimonious models.  
In particular, a library of candidate nonlinear terms $\boldsymbol{\Theta}(\mathbf{x})$ is constructed, and sparse regression is used to identify the few active terms in $\boldsymbol{\Theta}$ to approximate the function ${\bf f}$. %~\cite{Brunton2016pnas}.  

%\begin{figure}[tb]
%	\centering
%	\begin{overpic}[width=\textwidth]{Fig02}
%		\put(28.5,18.5){\scriptsize \cite{Brunton2017natcomm}}
%		\put(16,11.5){\scriptsize \cite{Kutz2016book}}
%		\put(20,5.5){\scriptsize \cite{Kaiser2017arxiv}}
%		\put(53,22.5){\scriptsize \cite{Loiseau2016jfm}}
%		\put(52,10.2){\scriptsize \cite{Mangan2016ieee}}
%		\put(79,20.5){\scriptsize \cite{Mangan2017prsa}}
%	\end{overpic}
%	\vspace{-.2in}
%	\caption{Schematic of the SINDYc algorithm (top row).  Active terms in a library of candidate nonlinearities are selected via sparse regression. 
%	Illustration of the modular nature of the SINDY (bottom row) with control framework and its ability to handle high-dimensional systems, limited measurements, known physical constraints, and model selection.}
%	\label{Fig:SINDYc:Schematic}	
%	\vspace{-.25in}
%\end{figure}

%\begin{figure*}[tb]
%	\centering
%	\includegraphics[width=\textwidth]{Fig02}
%	\vspace{-.2in}
%	\caption{Schematic of the SINDYc algorithm.  Active terms in a library of candidate nonlinearities are selected via sparse regression.  }
%	\label{Fig:SINDYc:Schematic}	
%\end{figure*}

SINDY with control (SINDYc) is based on the same assumption, that Eq.~\eqref{Eq:NonlinearDynamicsWithControl} only has a few active terms in the dynamics.  SINDY is readily generalized to include actuation, as this merely requires a larger library $ \boldsymbol{\Theta}(\mathbf{x},\mathbf{u})$ of candidate functions that include $\mathbf{u}$; these functions can include nonlinear cross terms in $\mathbf{x}$ and $\mathbf{u}$.  
 Thus, we measure $m$ snapshots of the state $\mathbf{x}$ and the input signal $\mathbf{u}$ in time and arrange these into two matrices:
%\begin{align}
%\mathbf{X}= \begin{bmatrix} \vline & \vline & & \vline \\
%\mathbf{x}_1 & \mathbf{x}_2 & \cdots & \mathbf{x}_m\\
%\vline & \vline & & \vline
%\end{bmatrix} ,\quad
%\mathbf{\Upsilon} = \begin{bmatrix} \vline & \vline & & \vline \\
%\mathbf{u}_1 & \mathbf{u}_2 & \cdots & \mathbf{u}_m\\
%\vline & \vline && \vline \end{bmatrix}.
%\end{align}
\begin{align}
\mathbf{X}= \begin{bmatrix} \mathbf{x}_1 & \mathbf{x}_2 & \cdots & \mathbf{x}_m
\end{bmatrix} ,\quad
\mathbf{U} = \begin{bmatrix} \mathbf{u}_1 & \mathbf{u}_2 & \cdots & \mathbf{u}_m
 \end{bmatrix}.
\end{align}
The library of candidate nonlinear functions $\boldsymbol{\Theta}$ may now be evaluated using the data $\mathbf{X}$ and $\mathbf{U}$:
%\begin{align}
%\mathbf{\Theta}^T(\mathbf{X},\boldsymbol{\Upsilon}) = \begin{bmatrix}
%\rule[2.2pt]{4em}{0.4pt} & \mathbf{1} & \rule[2.2pt]{4em}{0.4pt} \vspace{-.05in}\\ \vspace{-.05in}
%\rule[2.2pt]{4em}{0.4pt} & \mathbf{X} & \rule[2.2pt]{4em}{0.4pt}\\\vspace{-.05in}
%\rule[2.2pt]{4em}{0.4pt} & \boldsymbol{\Upsilon} & \rule[2.2pt]{4em}{0.4pt}\\\vspace{-.05in}
%%\rule[2.2pt]{4em}{0.4pt} & \mathbf{X}^2 & \rule[2.2pt]{4em}{0.4pt}\\
%\rule[2.2pt]{4em}{0.4pt} & \mathbf{X}\otimes\mathbf{X} & \rule[2.2pt]{4em}{0.4pt}\\\vspace{-.05in}
%\rule[2.2pt]{4em}{0.4pt} & \mathbf{X}\otimes\boldsymbol{\Upsilon} & \rule[2.2pt]{4em}{0.4pt}\\\vspace{-.05in}
%& \vdots &\\\vspace{-.05in}
%\rule[2.2pt]{4em}{0.4pt} & \sin(\mathbf{X}) & \rule[2.2pt]{4em}{0.4pt}\\\vspace{-.05in}
%\rule[2.2pt]{4em}{0.4pt} & \sin(\boldsymbol{\Upsilon}) & \rule[2.2pt]{4em}{0.4pt}\\\vspace{-.05in}
%\rule[2.2pt]{4em}{0.4pt} & \sin(\mathbf{X}\otimes\boldsymbol{\Upsilon}) & \rule[2.2pt]{4em}{0.4pt}\\\vspace{-.05in}
%%\rule[2.2pt]{4em}{0.4pt} & \sin(2\mathbf{X}) & \rule[2.2pt]{4em}{0.4pt}\\
%& \vdots &
%\end{bmatrix},
%\end{align}
%%% TRANSPOSE
\begin{align}
\hspace{-.1in}\boldsymbol{\Theta}(\mathbf{X},\mathbf{U}) = \begin{bmatrix}
\mathbf{1}^T \hspace{-.065in}& 
\hspace{-.065in}\mathbf{X}^T \hspace{-.065in}& 
\hspace{-.065in} \mathbf{U}^T \hspace{-.065in}& 
\hspace{-.065in}\left( \mathbf{X}\otimes\mathbf{X} \right)^T\hspace{-.065in}& 
\hspace{-.065in}\left( \mathbf{X}\otimes\mathbf{U} \right)^T\hspace{-.05in}
 \hspace{.015in} \cdots \hspace{.05in}
\hspace{-.025in}\sin(\mathbf{X})^T\hspace{-.05in}&
\hspace{-.05in}\sin(\mathbf{U})^T\hspace{-.05in}&
\hspace{-.05in} \sin(\mathbf{X}\otimes\mathbf{U})^T\hspace{-.025in} &
%\rule[2.2pt]{4em}{0.4pt} & \sin(2\mathbf{X}) & \rule[2.2pt]{4em}{0.4pt}\\
\hspace{-.1in} \cdots 
\end{bmatrix},
\end{align}
%\begin{align}
%\boldsymbol{\Theta}^T(\mathbf{X},\mathbf{U}) = \begin{bmatrix}
%\mathbf{1} \\ \vspace{-.05in}
%\mathbf{X}  \\\vspace{-.05in}
% \mathbf{U}\\\vspace{-.05in}
%%\rule[2.2pt]{4em}{0.4pt} & \mathbf{X}^2 & \rule[2.2pt]{4em}{0.4pt}\\
% \mathbf{X}\otimes\mathbf{X} \\\vspace{-.05in}
% \mathbf{X}\otimes\mathbf{U} \\\vspace{-.05in}
% \vdots \\\vspace{-.05in}
%\sin(\mathbf{X}) \\\vspace{-.05in}
%\sin(\mathbf{U}) \\\vspace{-.05in}
% \sin(\mathbf{X}\otimes\mathbf{U}) \\\vspace{-.05in}
%%\rule[2.2pt]{4em}{0.4pt} & \sin(2\mathbf{X}) & \rule[2.2pt]{4em}{0.4pt}\\
% \vdots 
%\end{bmatrix},
%\end{align}
where $\mathbf{x}\otimes\mathbf{y}$ defines the vector of all product combinations of the components in $\mathbf{x}$ and $\mathbf{u}$. 
\RefTwo{Although this definition includes repeated rows in $\boldsymbol{\Theta}$, in practice, the implementation is restricted to unique combinations.} 
A suitable library of candidate terms is crucial in the SINDYc algorithm.  
One strategy is to start with a basic choice, such as polynomials, and increase the complexity of the library by including other terms (trigonometric functions, etc.). %, which could also be systematically incorporated in a model selection approach.
It is also possible to incorporate partial knowledge of the physics (e.g.\ fluids vs. quantum mechanics) to decide on a library.

%It is possible to test different libraries (polynomials, trigonometric functions, etc.) and also incorporate partial knowledge of the physics (fluids vs. quantum mechanics, etc.).

The system in Eq.~\eqref{Eq:NonlinearDynamicsWithControl} can thus be written as:
\begin{align}
\dot{\mathbf{X}} = \boldsymbol{\Xi}\boldsymbol{\Theta}^T(\mathbf{X},\mathbf{U}).\label{Eq:SINDYc}
\end{align}
The time derivatives $\dot{\mathbf{X}}= \begin{bmatrix} \dot{\mathbf{x}}_1 & \dot{\mathbf{x}}_2 & \cdots & \dot{\mathbf{x}}_m \end{bmatrix}$, if not measured directly, are computed by numerical differentiation or approximated using the total variation regularized derivative~\cite{Rudin1992physd,Chartrand2011isrnam} if the data is noise-corrupted. 
The coefficients $\boldsymbol{\Xi}$ are \emph{sparse} for many dynamical systems.  
Therefore, we employ sparse regression to identify a sparse $\boldsymbol{\Xi}$ corresponding to the fewest nonlinearities in our library that give good model performance:  
\begin{align}
\boldsymbol{\xi}_k = \argmin_{\hat{\boldsymbol{\xi}}_k}\frac{\RefTwo{1}}{\RefTwo{2}}\|\dot{\mathbf{X}}_k-\hat{\boldsymbol{\xi}}_k\boldsymbol{\Theta}^T(\mathbf{X},\mathbf{U})\|_2^\RefTwo{2} + \lambda \|\hat{\boldsymbol{\xi}}_k\|_1,
\end{align}
where $\dot{\mathbf{X}}_k$ represents the $k$-th row of $\dot{\mathbf{X}}$ and $\boldsymbol{\xi}_k$ is the $k$-th row of $\boldsymbol{\Xi}$.  
%To approximate derivatives from noisy state measurements, the SINDY algorithm uses the total variation regularized derivative~\cite{Rudin1992physd,Chartrand2011isrnam}.

\begin{algorithm}
	\caption{\RefTwo{Sequentially thresholded least squares to learn the active library components.}}
	\label{Alg:SLST}
	\textbf{Input:} Time derivative $\dot{\mathbf{X}}$, library of candidate functions $\boldsymbol{\Theta}^T(\mathbf{X},\mathbf{U})$, thresholding parameter $\varepsilon$\\
	\textbf{Output:} Matrix of sparse coefficient vectors $\boldsymbol{\Xi}$\\\vspace{-.15in}
	\begin{algorithmic}[1]
		\Function{stls\_regression}{$\dot{\mathbf{X}}, \boldsymbol{\Theta}^T(\mathbf{X},\,\mathbf{U}), \varepsilon, N$}
		\State $\hat{\boldsymbol{\Xi}}^0 \leftarrow (\boldsymbol{\Theta}^T)^{\dagger}\dot{\mathbf{X}}$ \Comment{Initial least squares guess.}
		\While{not converged}
		\State $k \leftarrow k+1$
		\State ${\bf I}_{small} \leftarrow (\texttt{abs}(\hat{\boldsymbol{\Xi}})<\varepsilon)$ \Comment{Find small entries ...}
		\State $\hat{\boldsymbol{\Xi}}^k({\bf I}_{small}) \leftarrow0$ \Comment{... and threshold.}
		\For{all variables} 
		\State ${\bf I}_{big} \leftarrow\, \sim{\bf I}_{small}(:,ii)$ \Comment{Find big entries ...}
		\State $\hat{\boldsymbol{\Xi}}^k({\bf I}_{big},ii) \leftarrow 
		 (\boldsymbol{\Theta}^T(:,{\bf I}_{big}))^{\dagger}\dot{\mathbf{X}}(:,ii)$  \Comment{... and regress onto those terms.}
		\EndFor
		\EndWhile
		\EndFunction
	\end{algorithmic}
\end{algorithm}

The $\|\cdot\|_1$ term promotes sparsity in the coefficient vector $\boldsymbol{\xi}_k$.  This optimization may be solved using the LASSO~\cite{Tibshirani1996lasso} or the sequentially thresholded least squares procedure~\cite{Brunton2016pnas} \RefTwo{(see Alg.~\ref{Alg:SLST})}. 
\RefTwo{General conditions for the uniqueness of the $l_1$ relaxed solution have been provided in~\cite{tropp2006ieee}. In practice, these conditions may not be readily met, and false discoveries may occur, although they may be avoided under certain conditions~\cite{Su2016}. 
Specific conditions under which the sequentially thresholded least-squares algorithm in SINDy converges are provided in~\cite{Zhang2018arxiv}.	
More recently, convergence and recovery has been explored in a generalized framework for sparse relaxed regularized regression~\cite{Zheng2018arxiv}, for which SINDy constitutes a special case.
Conditions under which a model structure can be recovered from input--output data have also been examined in the context of identifiability~\cite{gevers2013cdc,alkhoury2017automatica}.}

%In the context of the $l_1$ relaxation problem, this has been theoretically examined by J. Tropp ({\it J. Tropp, ``Just relax'', ...}) outlining general conditions under which the problem is well-posed and 
%the sparse vector can be recovered. 
%Cand\`es et al. ({\it Cand\`s et al., ...}) demonstrated and analyzed under which conditions a unique solution can be found in practice.
%Specifically in the context of the SINDy architecture Schaefer et al. (...) provided theoretical convergence results ....
%More recently, Zheng et al. ({\it P. Zheng, T. Askham, S.~L. Brunton, J.~N. Kutz, A.~Y. Aravkin, ``A Unified Framework for Sparse Relaxed Regularized Regression: SR3'', arXiv:1807.05411}) 

The parameter $\lambda$ (or equivalently $\varepsilon$ in Alg.~\ref{Alg:SLST}) is selected to identify the Pareto optimal model that best balances model complexity with accuracy.  A coarse sweep of $\lambda$ is performed to identify the rough order of magnitude where terms are eliminated and where error begins to increase.  Then this parameter sweep may be refined, and the models on the Pareto front are evaluated using information criteria~\cite{Mangan2017prsa}.  
It is interesting to note, that a similar idea, identifying active components in ${\bf f}$ from a library of candidate functions using sparse regularization, was discarded in favor of a Bayesian formulation, as the non-orthogonality of the columns in the library was seen as problematic~\cite{pan2016ieee}. However, as in~\cite{Brunton2016pnas}, we will demonstrate here the effectiveness of the approach.

Since the original SINDY paper~\cite{Brunton2016pnas}, it has been extended to include constraints and known physics~\cite{Loiseau2016jfm}, for example to enforce energy preserving constraints in an incompressible fluid flow.  
SINDY has also been extended to high-dimensional systems, by identifying dynamics on principal components~\cite{Brunton2016pnas}, learning partial differential equations~\cite{Rudy2017sciadv,Schaeffer2017prsa}, and extracting dynamics 
%from a few limited measurements~\cite{Loiseau2017arxiv} or 
on delay coordinates~\cite{Brunton2017natcomm}. 
Robust variants of SINDY have been formulated to identify models despite large outliers and noise~\cite{Tran2016arxiv,Schaeffer2017pre}.  

\subsubsection{Discovering discrete-time dynamics}
In the original SINDY algorithm, it was shown that it is possible to identify discrete-time models of the form 
%
%\begin{align}
$\bx_{k+1} = \bF(\bx_k)$.
%\end{align}
%
It is also possible to extend SINDY to identify discrete-time models with inputs and control:
\begin{align}\label{Eq:NLDiscrete}
\bx_{k+1} = \bF(\bx_k,\bu_k).
\end{align}
Instead of computing derivatives, we collect a matrix $\bX'$ with the columns of $\bX$ advanced one timestep:
%%
%\begin{align}
%\bX' =\begin{bmatrix} \vline & \vline & & \vline \\
%\mathbf{x}_2 & \mathbf{x}_3 & \cdots & \mathbf{x}_{m+1}\\
%\vline & \vline & & \vline
%\end{bmatrix}.
%\end{align}
%%
%
%\begin{align}
$\bX' =\begin{bmatrix} \mathbf{x}_2 & \mathbf{x}_3 & \cdots & \mathbf{x}_{m+1}
\end{bmatrix}$.
%\end{align}
%
Then, the dynamics may be written as
\begin{align}
\bX' = \boldsymbol{\Xi}\boldsymbol{\Theta}^T(\mathbf{X},\mathbf{U}),
\end{align}
and the regression problem becomes 
\begin{align}
\boldsymbol{\xi}_k = \argmin_{\hat{\boldsymbol{\xi}}_k}\frac{\RefTwo{1}}{\RefTwo{2}}\|{\mathbf{X}'}_k-\hat{\boldsymbol{\xi}}_k\boldsymbol{\Theta}^T(\mathbf{X},\mathbf{U})\|_2^\RefTwo{2} + \lambda \|\hat{\boldsymbol{\xi}}_k\|_1.
\end{align}

%The discrete-time dynamics in Eq.~\eqref{Eq:NLDiscrete} are more general than continuous-time dynamics in Eq.~\eqref{Eq:NonlinearDynamicsWithControl} and don't require the computation of derivatives. 
%%
%\begin{align}
%
%\end{align}
%%

%\begin{itemize}
%	\item sparsity-promoting l1, parameter sweep
%	\item Choosing library candidate functions
%	\item incorporate partial knowledge of physics
%	\item partial measurements, high-dimensional systems (SVD/POD, feature engineering)
%	\item Model selection
%	\item Convergence to DMDc
%	\item derivatives from noisy measurements using total variation regularized derivative
%\end{itemize}

\subsubsection{Relationship to dynamic mode decomposition}
The SINDY regression is related to the dynamic mode decomposition (DMD), which originated in the fluids community to extract spatiotemporal coherent structures from large fluid data sets~\cite{Rowley2009jfm,schmid:2010,Tu2014jcd,Kutz2016book}.    
DMD modes are spatially coherent and oscillate at a fixed frequency and/or growth or decay rate.  
Since fluids data is typically high-dimensional, DMD is built on the proper orthogonal decomposition (POD)~\cite{HLBR_turb}, effectively recombining POD modes in a linear combination to enforce the temporal coherence.  
The dynamic mode decomposition has been applied to a wide range of problems including fluid mechanics, epidemiology, neuroscience, robotics, finance, and video processing~\cite{Kutz2016book}.  
Many of these applications have the ultimate goal of closed-loop feedback control.

In DMD, a similar regression is performed to identify a linear discrete-time model $\mathbf{A}$ mapping $\bX$ to $\bX'$:
\begin{align}
\mathbf{X}' = \mathbf{A} \mathbf{X}.
\end{align}
Thus, SINDY reduces to DMD if formulated in discrete-time, with linear library elements in $\mathbf{\Theta}$, and without a sparsity-promoting $L_1$ penalty term, i.e.\ $\lambda=0$.  

DMD was recently extended to include actuation inputs by Proctor et al~\cite{Proctor2016siads}, to disambiguate the effect of internal dynamics and control.  
In dynamic mode decomposition with control (DMDc), a similar regression is formed, but with the actuation input matrix $\mathbf{U}$:
\begin{align}
\mathbf{X}' = \mathbf{A} \mathbf{X} + \mathbf{B}\mathbf{U}.
\end{align}
Thus, SINDY with control similarly reduces to DMDc under certain conditions.  
In this work, we will use DMDc and SINDYc to discover dynamics for model predictive control.  
The DMDc algorithm has also been shown to be related to other subspace identification methods, such as the eigensystem realization algorithm~\cite{ERA:1985}, but designed for high-dimensional input--output data. 

It is interesting to note that the extended DMD (eDMD)~\cite{Williams2015jnls} regression is performed on the nonlinear library $\mathbf{\Theta}(\mathbf{X}') = \mathbf{A}\mathbf{\Theta}(\mathbf{X})$, and an $L_1$ penalty may also be added.  
Extended DMD may also be modified to incorporate actuation inputs, and these models have recently been used effectively for model predictive control~\cite{korda2016_a}.

\subsubsection{Identification of dynamics with feedback control}
If the input $\mathbf{u}$ corresponds to feedback control, so that $\mathbf{u} = \mathbf{K}(\mathbf{x})$, then it is impossible to disambiguate the effect of the feedback control $\mathbf{u}$ with internal feedback terms $\mathbf{K}(\mathbf{x})$ within the dynamical system; namely, the SINDYc regression becomes ill-conditioned.  In this case, we may identify the actuation $\mathbf{u}$ as a function of the state:
\begin{align}
\mathbf{U} = \boldsymbol{\Xi}_u\boldsymbol{\Theta}^T(\mathbf{X}).
\end{align}  
To identify the coefficients $\boldsymbol{\Xi}$ in Eq.~\eqref{Eq:SINDYc}, we perturb the signal $\mathbf{u}$ to allow it to be distinguished from $\mathbf{K}(\mathbf{x})$ terms.  
This may be done by injecting a sufficiently large white noise signal, or occasionally kicking the system with a large impulse or step in $\mathbf{u}$.  
An interesting future direction would be to design input signals that \emph{aid} in the identification of the dynamical system in Eq.~\eqref{Eq:NonlinearDynamicsWithControl} by perturbing the system in directions that yield high-value information.

\begin{figure}[tb]
	\centering
	\begin{overpic}[width=\textwidth]{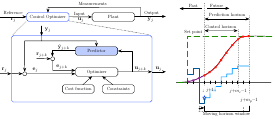}
%		\put(0,0){(a)}
%		\put(60,0){(b)}
		\put(0,2){(a)}
		\put(60,2){(b)}
	\end{overpic}	
	\vspace{-.3in}
	\caption{
	Schematics for (a) the control loop and (b) the receding horizon of MPC.
	Full state measurements $\by = \bx$, and $\hat\by = \hat\bx$ are considered as the output in the examples. 
	Starting from the most recent measurement, the control input sequence (light blue solid) is optimized over the control horizon based on predicted future outputs (red solid) to drive the system to the set point (green dashed). 
	The first input (blue star) in the sequence is enacted.	
%	Past measurements and control inputs are displayed as purple and blue solid lines, respectively.
%	Schematics for model predictive control:  (a) control loop, and (b) receding horizon. 
%	In the examples demonstrating SINDY-MPC, the output is the full state $\by = \bx$, and $\hat\by = \hat\bx$ is the model prediction. Starting from the most recent measurement, the control input sequence (light blue solid) is optimized over the control horizon based on predicted future outputs (red solid) in order to drive the system to the set point (green dashed). 
%	The first input (blue star) in the sequence is enacted.	
%	Past measurements and control inputs are displayed as purple and blue solid lines, respectively.
	}
	\label{Fig:MPC:MPC_ControlLoop}	
	\vspace{-.2in}
\end{figure}

%\begin{figure*}[tb]
%	\centering
%	\includegraphics[width=0.75\textwidth]{Fig03}
%	\vspace{-.1in}
%	\caption{Schematic for model predictive control.  In SINDY-MPC, the output is the full state $\by = \bx$, and $\hat\by = \hat\bx$ is the model prediction. }
%	\label{Fig:MPC:MPC_ControlLoop}	
%\end{figure*}

%\begin{figure}[tb]
%	\centering
%	\includegraphics[width=0.375\textwidth]{Fig04}
%	\vspace{-.2in}
%%	\caption{Receding horizon control: The control input sequence (light blue solid) is optimized over the control horizon based on predicted future outputs (red solid line with dots) in order to drive the system to the set point (green dashed). 
%%		The first input (blue/gray star) in the sequence is applied in the plant.	
%%		Past measurements and control inputs are displayed as purple and blue solid lines, respectively.}
%	\caption{Receding horizon control: The control input sequence (light blue solid) is optimized over the control horizon based on predicted future outputs (red solid) in order to drive the system to the set point (green dashed). 
%		The first input (blue star) in the sequence is enacted.	
%		Past measurements and control inputs are displayed as purple and blue solid lines, respectively.}
%	\label{Fig:MPC:RecedingHorizon}	
%\end{figure}

%\cleardoublepage
\subsection{Model predictive control}
In this section, we outline the control problem and summarize key results in MPC, which is shown schematically in Fig.~\ref{Fig:MPC:MPC_ControlLoop}.   
Model predictive control solves an optimal control problem over a receding horizon, subject to system constraints, to determine the next control action.  
This optimization is repeated at each new timestep, and the control law is updated, as shown in Fig.~\ref{Fig:MPC:MPC_ControlLoop}(b).

The receding horizon control problem can generally be formulated as an open-loop optimization at each step, which determines the optimal sequence of control inputs $\bu(\cdot\vert\bx_j):= \{\bu_{j+1},\ldots,\bu_{j+k},\ldots,\bu_{j+m_c}\}$ over the control horizon $T_c=m_c\Delta t$ given the current measurement $\bx_j$  that minimizes a cost $J$ over the prediction horizon $T_p=m_p\Delta t$; $\Delta t$ is the timestep of the model, which may be different from the sampling time of measurements.  
The control horizon is generally less than or equal to the prediction horizon, so that $T_c \leq T_p$; if $T_c<T_p$, then the input $\bu$ is assumed constant thereafter.
\RefTwo{The first control value $\bu_{j+1}$ %$\bu_{j+1}:=\bu(\bu_{j+1}\vert \bx_j)$
is then applied}, and the optimization is reinitialized and repeated at each subsequent timestep
\RefTwo{to solve for the unknown sequence $\bu(\cdot\vert\bx_j)$}.  
This results in an 
implicit feedback control law
\RefTwo{
\begin{align}
\bK(\bx_j) = \bu(j+1\vert\bx_j) = \bu_{j+1},
\end{align}
}
where $\bu_{j+1}$ is the first in the optimized actuation sequence starting at the initial condition $\bx_j$.  

The cost optimization at each timestep is given by:
\begin{align}
\hspace{-.1in}\min\limits_{\hat{\bu}(\cdot\vert\bx_j)} J(\bx_j) =& \min\limits_{\hat{\bu}(\cdot\vert\bx_j)} \left[ \vert\vert \hat\bx_{j+m_p} - \bx_{m_p}^{*}\vert\vert_{\bQ_{m_p}}^2 + \sum\limits_{k=0}^{m_p-1}\,  \vert\vert \hat\bx_{j+k}-\bx_{k}^{*} \vert\vert_{\bQ}^2  \right. \notag\\ 
&+ \left.\sum\limits_{k=1}^{m_c-1}\, \left(\vert\vert \hat\bu_{j+k} \vert\vert_{\bR_u}^2 + \vert\vert \Delta\hat\bu_{j+k} \vert\vert_{\bR_{\Delta u}}^2 \right)\right],\label{Eqn:MPC:OptimalControlProblem}
\end{align}
subject to the discrete-time system dynamics with $\hat{\bf F}:\mathbb{R}^{n}\times \mathbb{R}^q\rightarrow\mathbb{R}^n$
\begin{equation}\label{Eq:NonlinearDynamicsWithControl_Model}
\hat{\bx}_{k+1} =\hat{{\bf F}}(\hat{\bx}_k,\bu_k),
\end{equation}
the input constraints,
\begin{subequations}
	\begin{align}
	\Delta \bu_{min} \leq \Delta &\bu_k \leq \Delta \bu_{max},\\
	\bu_{min} \leq  &\bu_k \leq \bu_{max},
	\end{align}
\end{subequations}
and possibly additional equality or inequality constraints on the state and input.  
Here, we assume the availability of full-state measurements $\by = \bx$.
The cost functional $J$ penalizes deviations of the predicted state $\hat\bx_k$ along the trajectory $\bx_{k}^{*}$ and also includes a terminal cost at $\hat\bx_{m_p}$.  
Expenditures of the input $\bu_k$ and input rate $\Delta \bu_k = \bu_k-\bu_{k-1}$ are also penalized.  
Each term is computed as the weighted norm of a vector, i.e., $\vert\vert \bx \vert\vert_{\bQ}^2 := \bx^T\bQ\bx$.  
The weight matrices $\bQ\geq 0$, $\bQ_{m_p}\geq 0$, $\bR_u>0$ and $\bR_{\Delta u}>0$ are positive definite and positive semi-definite, respectively. 
Note that the model prediction $\hat\bx_k$, which is forecast, may differ from the true measured state $\bx_k$.  

\RefTwo{The dynamics are given by the identified SINDYc model, e.g. $\dot{\bf x} = {\bf F}({\bf x},{\bf u}) = \boldsymbol{\Xi}\boldsymbol{\Theta}^T({\bf x},{\bf u})$; $\hat{{\bf F}}$ represents a discrete-time or discretized SINDYc model. While the model and the control law may be learned simultaneously, we adopt a two-stage process, where the model is first learned from data and then used in the control optimization with MPC. 
A joint optimization of the model and the control law may be challenging, as the particular control action depends on the model. 
However, it may be possible to develop a streaming algorithm to adapt the model to abrupt system changes~\cite{Quade2018arxiv}, iterating between model identification and control optimization.  }

MPC is one of the most powerful model-based control techniques due to the flexibility in the formulation of the objective functional, the ability to add constraints, and extensions to nonlinear systems. 
%The most challenging aspects of MPC involve 
%the fast and robust real-time optimization and
The most challenging aspect of MPC involves
the identification of a dynamical model that accurately and efficiently represents the system behavior when control is applied. 
% In this work, we focus on the latter.
If the model is linear, minimization of a quadratic cost functional subject to linear constraints results in a tractable convex problem.  
Nonlinear models may yield significant improvements; however, they render MPC a nonlinear program, which can be expensive to solve, making it particularly challenging for real-time control. 
%However, nonlinear models render MPC a nonlinear program, which can be expensive to solve, making it particularly challenging for real time control. 
%Improvements from nonlinear models may be small for processes at operating conditions with small disturbances.
%However, improved forecasting using nonlinear models becomes key when processes are operated at various operating conditions over short time scales or over large state space regions,  when large disturbances occur, 
\RefTwo{Conditions on the well-posedness of the problem and existence and uniqueness of the solution of the nonlinear optimization problem are, e.g., provided in~\cite{johansen2011intro}.}
Fortunately, improvements in computing power and advanced algorithms are increasingly enabling nonlinear MPC for real-time applications.

%% file: Sec41.tex
\begin{figure}[tb]
	\centering
	\includegraphics[width=\textwidth]{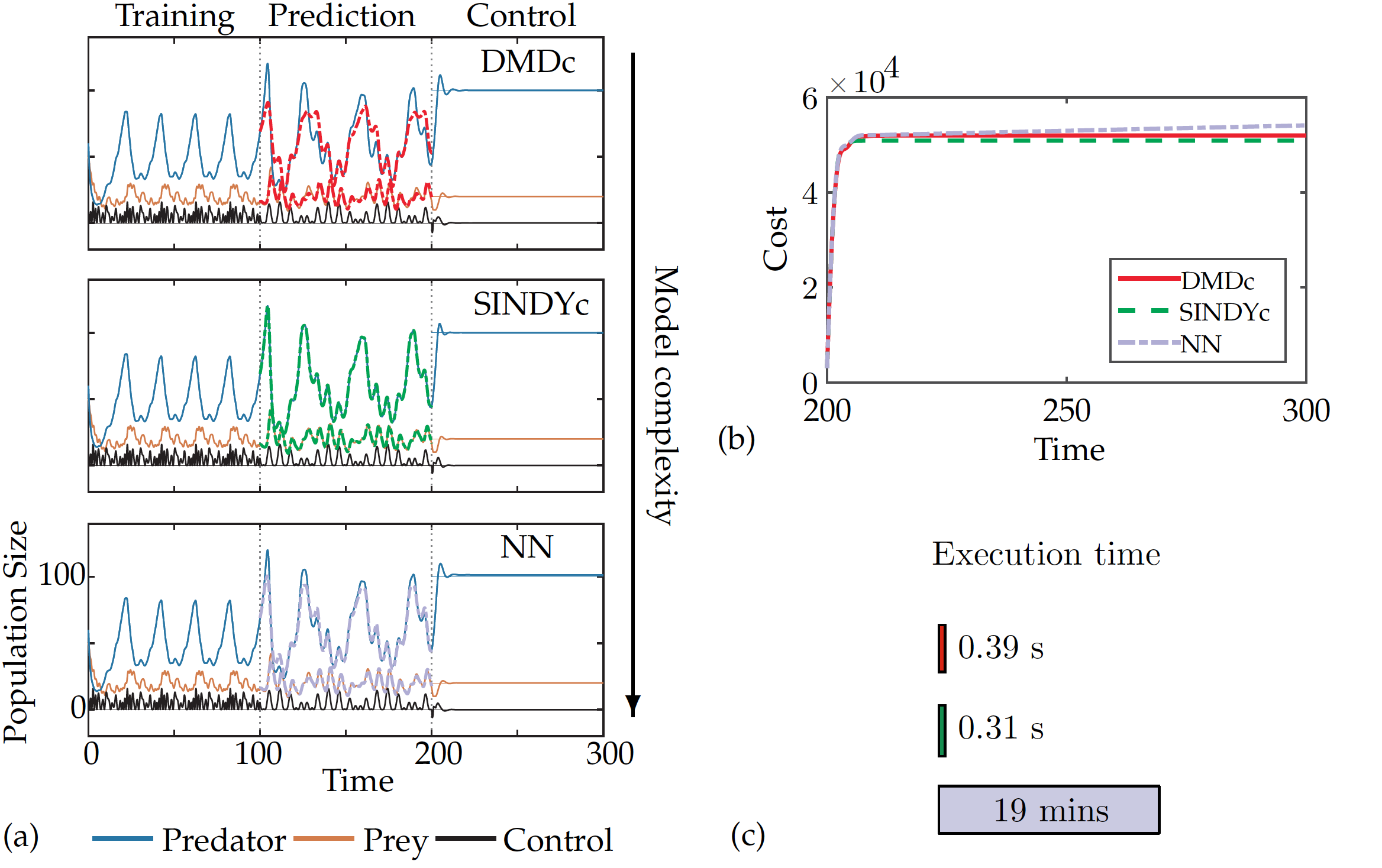}
	\vspace{-.2in}
	\caption{Prediction and control performance for the Lotka-Volterra system:
		(a) time series of states and input during training, validation, and control stage, and
		(b) cumulative cost, and 
		(c) execution time of the MPC optimization procedure.}
	\label{Fig:LOTKA:MPC_Comparison}
	\vspace{-.3in}
\end{figure}

\section{A simple model for population dynamics}\label{Sec:Lotka-Volterra}
We first demonstrate the SINDY-MPC architecture on the Lotka-Volterra system, 
a two-dimensional, weakly nonlinear dynamical system, describing the interaction between two competing populations.  
These dynamics may represent two species in biological systems, 
%the dynamic 
competition in stock markets~\cite{Lee2005tfsc}, 
and %they 
can be modified to study the spread of infectious diseases~\cite{Venturino1994jm}. 
We will consider more sophisticated examples in the following sections.

The dynamics of the prey and predator populations, $x_1$ and $x_2$, respectively, are given by
\begin{subequations}\label{Eqn:Lotka-Volterra}
	\begin{align}
	\dot{x}_1 &= ax_1 - bx_1 x_2\\
	\dot{x}_2 &= -cx_2 + d x_1 x_2  + u,
	\end{align}
\end{subequations}
where the constant parameters $a = 0.5$, $b=0.025$, $c=0.5$, and $d=0.005$ represent the growth/death
rates, the effect of predation on the prey population, and the growth of predators based on the size of the prey population.
%Here, the actuation input $u$ can assume negative values representing a constant-yield harvesting or positive values depicting the release of new predators into the wild.
% control $u$ represents a proportional removal of predator population ({\it http://www.scielo.br/pdf/ca/v16n2/a02v16n2.pdf})
% in 2nd eqn:  - x_2 u
The unforced system exhibits a limit cycle behavior, where the predator lags the prey, and a critical point $\bx^{crit} =(g/d\;\; a/b)^T$, where the population sizes of both species are in balance.
The control objective is to stabilize this fixed point.  
Here, the timestep $\Delta t = 0.1$ of the system and the model are equal,  
the weight matrices are $\bQ = (\begin{smallmatrix} 1 & 0\\ 0 & 1 \end{smallmatrix})$ and $R_u = R_{\Delta u} = 0.5$, and the actuation input is limited to $u\in[-20,20]$.
The control and prediction horizons are $m_p = m_c = 10$ unless otherwise noted.
%The possibility to add constraints is one reason why MPC has become so successful.
We apply an additional constraint on $u$, so that $x_2$ does not decrease below $10$, to enforce a minimum population size required for recovery. 

\begin{figure}[tb]
	\centering
	\includegraphics[width=\textwidth]{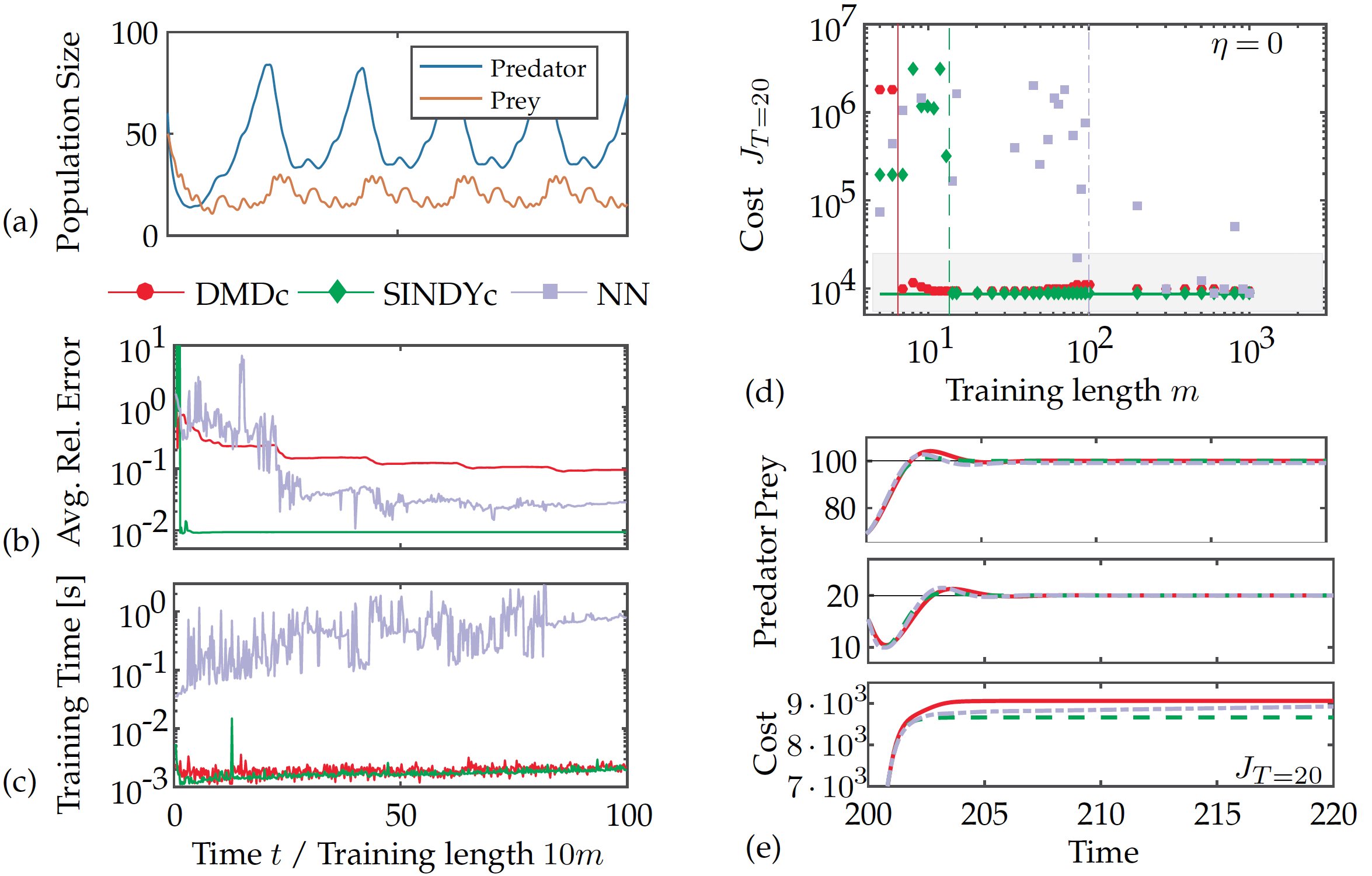}
	\vspace{-.2in}
	\caption{Crossvalidated prediction error and control performance for increasing length of training data: 
		(a) training time series, % of the training data, 
		%		(b) mean squared error, 
		(b) average relative prediction error, 
		(c) model training time in seconds,
		(d) terminal cumulative cost over 20 time units, and % for a representative set of models, and
		(e) time series of states and cost of the best model for each model type ($m^{DMDc} = 20$, $m^{SINDYc} = 85$, $m^{NN} = 10^3$). 
		From $m= 14$ onwards, SINDYc yields highly performing models in MPC, outperforming all other models. 
		Outside the shaded region, models perform significantly worse or even diverge.
	}
	\label{Fig:LOTKA:IncreasingTrainingLength}
	\vspace{-.3in}
\end{figure}

To assess the performance and capabilities of the SINDY-MPC architecture, SINDYc is compared with two representative data-driven models:
dynamic mode decomposition with control (DMDc) and a multilayer neural network (NN), which can represent any continuous function under mild conditions~\cite{hornik1989multilayer}.  
The results are displayed in Fig.~\ref{Fig:LOTKA:MPC_Comparison}. 
The first $100$ time units are used to train the models with a phase-shifted sum of sinusoids as input, a so-called Schroeder sweep~\cite{schroeder1970ieee}, after which the predictive capabilities of these models are validated using sinusoidal forcing with $u(t) = (2 \sin(t) \sin(t/10))^2$ on the next $100$ time units. 
Different actuation inputs are used during the training and validation stages to assess the models' ability to generalize.
Thereafter, MPC is applied for $100$ time units using a prediction and control horizon of $m_p=m_c=5$.
SINDYc shows the best prediction and control performance, followed by DMDc and the neural network (due to its steady-state error).  
The neural network has $1$ hidden layer with $10$ neurons, which is the best trade-off between model complexity and accuracy; increasing the number of neurons or layers has little impact on the prediction performance.
Further, hyperbolic tangent sigmoid activation functions are employed.
It is first trained as a feedforward network using the Levenberg-Marquardt algorithm and then closed.
If the data is corrupted by noise, a Bayesian regularization is employed, which requires more training time but improves robustness.
%Following heuristic approaches in~\cite{???}, we choose $1$ hidden layer with $10$ neurons in the neural network, which represents the best trade-off between model complexity and accuracy.
%Increasing the number of neurons or layers had hardly any impact on improving the prediction performance.
While the neural network exhibits a similar control performance, the execution time of SINDYc is $37$ times faster, which is particularly critical in real-time applications.
For a fair comparison, all methods are compared using the same optimization routine based on interior-point methods via Matlab's \texttt{fmincon}. Thus, it would be possible to reduce the time for the linear system further.

\begin{figure}[tb]
	\centering
	\includegraphics[width=\textwidth]{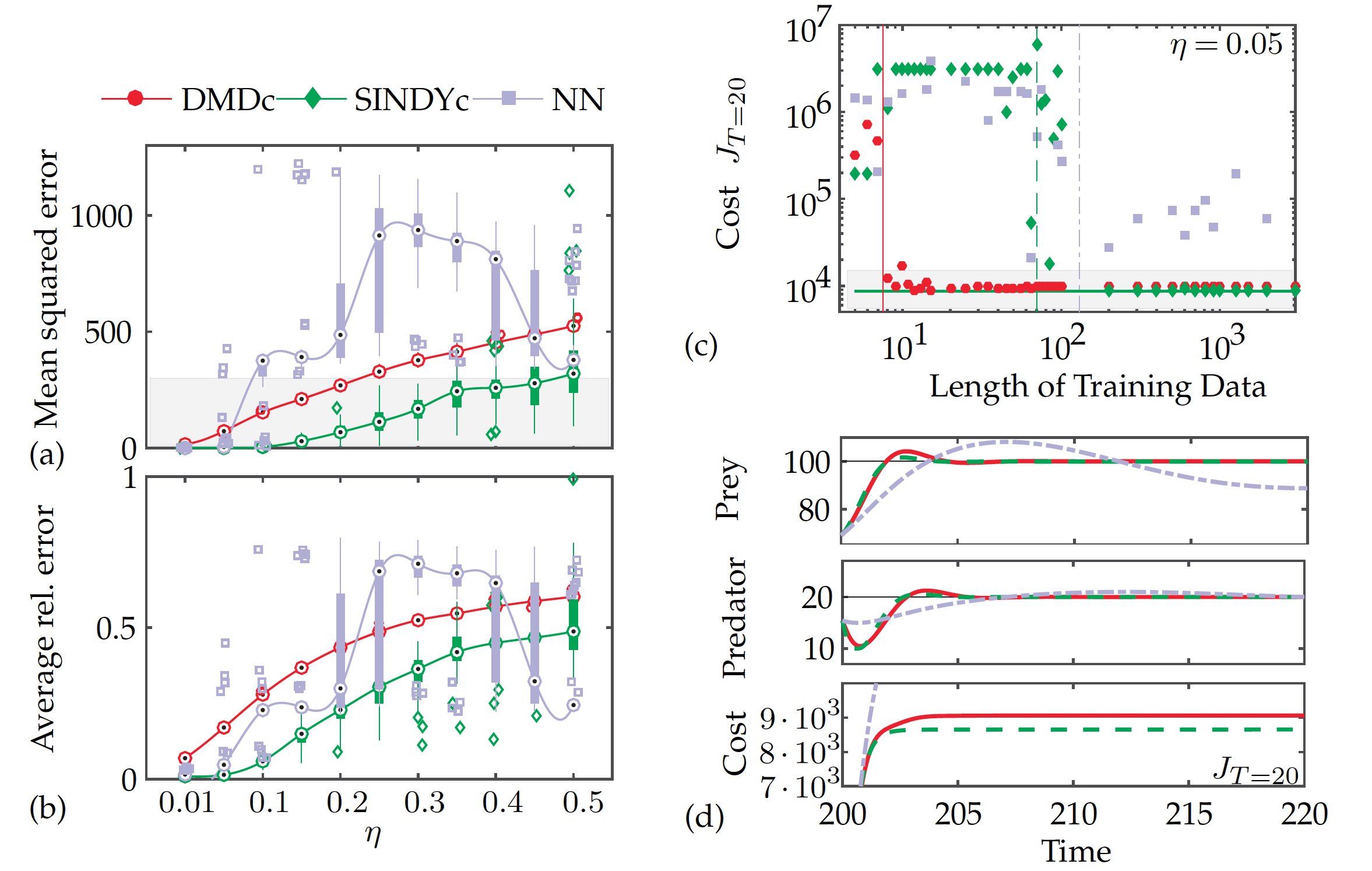}
	\vspace{-.3in}
	\caption{
		Crossvalidated prediction and control performance for increasing measurement noise:
		%using a Schroeder-type excitation signal: 
		(a) mean squared and (b) average relative prediction error,  
		(c) terminal cumulative cost over 20 time units, and
		(d) time series of states and cost of the best models
		($m^{DMDc} = 12$, $m^{SINDYc} = 1250$, $m^{NN} = 65$). 
		The control performance is shown for increasing length of noise-corrupted training data with $\eta = 0.5$.
		From $m = 200$ onwards, SINDYc yields highly performing models, outperforming all other models. 
		Statistics are shown for 50 noise realizations each. 
		Above the shaded region, models in most realizations do not have any predictive power.
		%(except for one model that performs slightly worse than the best DMDc model).
		}
	\label{Fig:LOTKA:DependencyOnTrainingLength_WithNoise}
\vspace{-.3in}	
\end{figure}

In practice, measurements are generally affected by noise.
We examine the robustness of these models for increasing noise corruption of the state measurements, 
i.e.\ ${\bf y} = {\bf x} + {\bf n}$ where ${\bf n}\in\mathcal{N}(0,{\sigma}^2)$ with standard deviation $\sigma$.
Crossvalidated prediction performance for different noise magnitudes $\eta=\sigma/\max(\mathrm{std}(x_i))\in(0.01,0.5)$, where $\mathrm{std}$ denotes standard deviation,
is displayed in Fig.~\ref{Fig:LOTKA:DependencyOnTrainingLength_WithNoise}(a-b).
As expected, the performance of all models decreases with increasing noise magnitude.
SINDYc generally outperforms DMDc and neural network models, exhibiting a slower decline in performance for low and moderate noise levels. 
Sparse regression is known to improve robustness to noise and prevent overfitting.  
The large fluctuation in the neural network performance are due to its strong dependency on the initial network weights.

The amount of data required to train an accurate model is particularly crucial in real-time applications, where abrupt changes or actuation may render the model invalid and rapid model updates are necessary.
Figure~\ref{Fig:LOTKA:IncreasingTrainingLength}(a-c) shows the average relative prediction error on $100$ time units used for validation, and the training time for increasing lengths of training data. 
The effect of the training length on the control performance (evaluated over 20 time units) is shown in Fig.~\ref{Fig:LOTKA:IncreasingTrainingLength}(d-e).
For small amounts of data, the sparsity-promoting parameter $\lambda$ in SINDYc is reduced by a factor of $10$ until a non-zero entry appears.
In the low-data limit, a highly predictive SINDYc model can be learned, discovering the true governing equations within machine precision. 
Significantly larger amounts of data are required to train an accurate neural network model, although with enough data it outperforms DMDc.  
DMDc models may be useful in the extremely low-data limit, before enough data is available to characterize a SINDYc model. 
The training times of SINDYc and DMDc models increase slightly with the amount of data, but they require about two orders of magnitude less time than neural network models. 
SINDYc's intrinsic robustness to overfitting renders all models from $m_{train}=14$ on as having the best control performance compared with the overall best performing DMDc and neural network models. 
In contrast, DMDc shows a slight decrease in performance due to overfitting
and the neural network's dependency on the initial network weights detrimentally affects its performance. 
It is interesting to note that the control performance is generally less sensitive than the long-term prediction performance shown in Fig.~\ref{Fig:LOTKA:IncreasingTrainingLength}(b-c).
Even a model with moderately low predictive accuracy may perform well in MPC.

In Fig.~\ref{Fig:LOTKA:DependencyOnTrainingLength_WithNoise}(c-d) we show the same analysis but with noise-corrupted training data. 
We assume no noise corruption during the control stage.
For each training length, the best model out of $50$ noise realizations is tested for control.
DMDc and SINDYc models both require slightly more data to achieve a similar performance as without noise.
Note that neural network models perform significantly worse when trained on noise-corrupted data.

%\begin{itemize}
%%	\item DMDc assuming unknown ${\bf B}$
%%	\begin{itemize}
%%		\item New state ${\bf y} = {\bf x} - {\bf x}^{ref}$ 
%%	\end{itemize}
%%	\item SINDy with  penalization $\lambda = 0.001$
%	\item NARX
%	\begin{itemize}
%		\item 1 hidden layers with 10 neurons (heuristic rule of thumb = at least 2x size of input)
%		\item Up to 1000 iterations in feedforward training
%		\item random initial weights (is generally a bad idea, improve? --> no)
%		\item hyperbolic tangent sigmoid activation function
%	\end{itemize}
%\end{itemize}

%% file: Sec42.tex
\begin{figure}[tb]
	\centering
	\includegraphics[width=\textwidth]{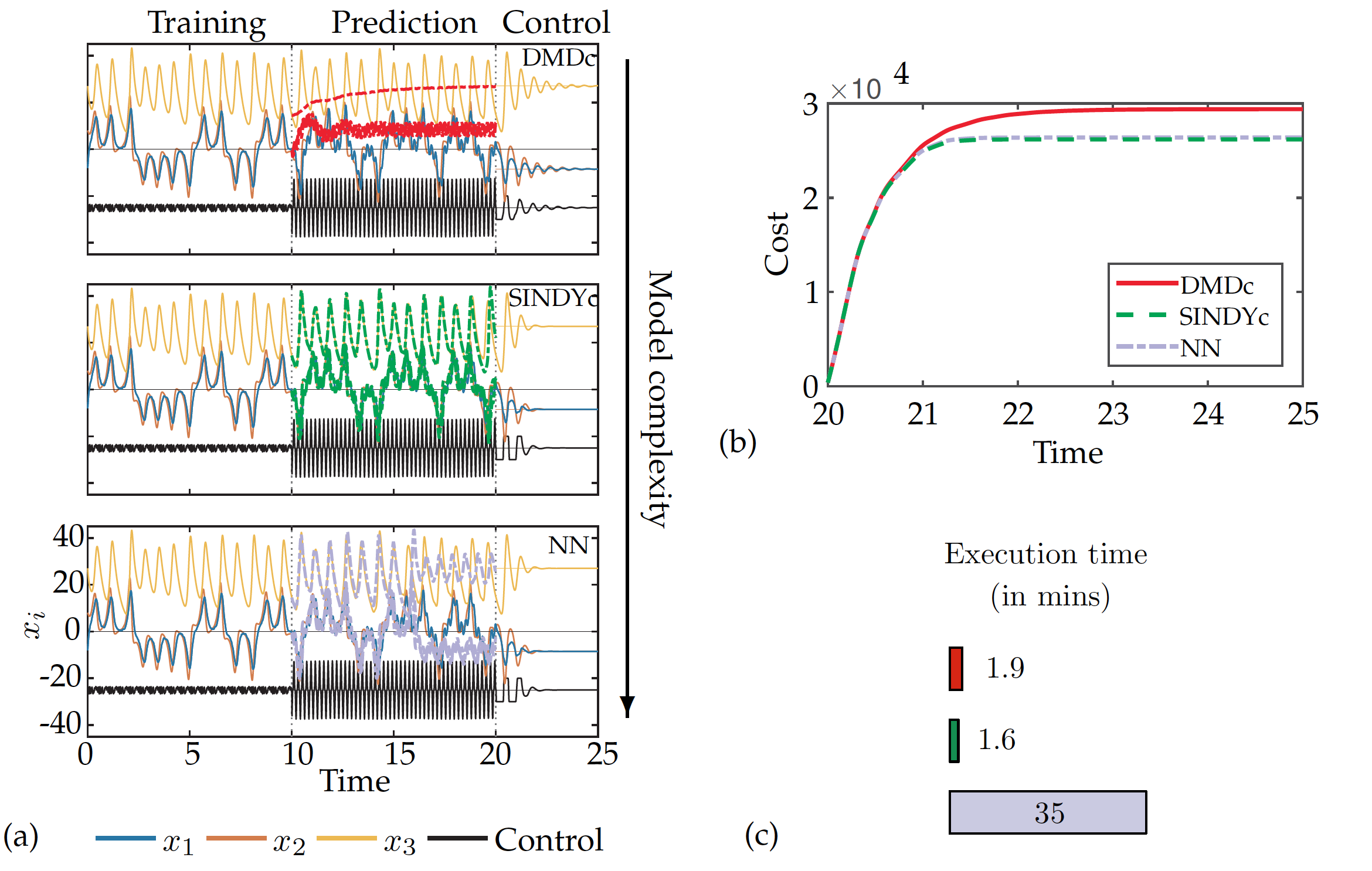}
	\vspace{-.3in}
	\caption{Prediction and control performance for the chaotic Lorenz system:
		(a) time series of the states and input (shifted to $-25$ and scaled by $10$ to improve readability) during training, validation, and control stage, and
		(b) cumulative cost, and 
		(c) execution time of the MPC optimization.}
	\label{Fig:LOREZ:MPC_Comparison}
	\vspace{-.3in}
\end{figure}

\section{Chaotic Lorenz system}\label{Sec:Lorenz}
In this section, we demonstrate the SINDY-MPC architecture on the chaotic Lorenz system, 
a prototypical example of chaos in dynamical systems.
The Lorenz system represents the Rayleigh-B\'enard convection in fluid dynamics as proposed by Lorenz~\cite{Lorenz1963jas}, 
but has also been associated with lasers, dynamos, and chemical reaction systems. 
The Lorenz dynamics are given by
\begin{subequations}\label{Eqn:Lorenz}
	\begin{align}
	\dot{x}_1 &= \sigma(x_2-x_1)  + u\\
	\dot{x}_2 &= x_1(\rho - x_3) - x_2 \\
	\dot{x}_3 &= x_1 x_2 - \beta x_3 
	\end{align}
\end{subequations}
with system parameters $\sigma = 10$, $\beta = 8/3$, $\rho = 28$, and control input $u$ affecting only the first state.
%For these values of the parameters, the attractor possesses three unstable fixed points at
%$(0,0,0)^T$ and $(\pm \sqrt{72}, \pm\sqrt{72},27)^T$. 
A typical trajectory oscillates alternately around the two weakly unstable fixed points $(\pm \sqrt{72}, \pm\sqrt{72},27)^T$. % and exhibits a strong sensitivity to initial conditions.
The chaotic motion of the system implies a strong sensitivity to initial conditions, i.e.\ small uncertainties in the state will grow exponentially with time. 
This represents a particularly challenging problem for model identification and subsequent control, as measurement and model uncertainty both lead to long-time forecast error. %
%Discovering the model from noisy measurements renders this issue more severe as model uncertainties also affect the forecast accuracy.
 
The control objective is to stabilize one of these fixed points.
In general, the timestep of the model is chosen to balance the control horizon, the length of the sequence of control inputs to be optimized, and prediction accuracy.
%, while assuring the model is as predictive as possible for the given timestep. 
Here, the system timestep  is $\Delta t^{sys} = 0.001$ and the model timestep is $\Delta t^{model} = 0.01$. The control input is determined every $10$ system timesteps and then held constant.
%In this example, the chosen timesteps serve well these purposes.
The weight matrices are $\bQ = \bI_{3}$, where $\bI_n$ denotes a $n\times n$ identity matrix, 
%(\begin{smallmatrix} 1 & 0 & 0\\ 0 & 1 & 0\\ 0 & 0 & 1\end{smallmatrix})$
$R_u = R_{\Delta u} = 0.001$, and the actuation input is limited to $u\in[-50,50]$.
The control and prediction horizon is $m_p = m_c = 10$ and the sparsity-promoting parameter in SINDYc is $\lambda = 0.1$, unless otherwise noted. 
For all cases we assume access to full-state information.

We compare the prediction and control performance of the SINDYc model with DMDc and neural network models. 
DMDc is trained to model the deviation from the goal state by constructing the regression model based on data from which the goal state has been subtracted. 
A less na\"ive approach would partition the trajectory into two bins, e.g.\ based on negative and positive values of $x_1$, and estimate two models for each goal state separately.
The neural network consists of $1$ hidden layer with $10$ neurons and employs hyperbolic tangent sigmoid activation functions.
Cross-validated prediction and control performance for the Lorenz system are displayed in Fig.~\ref{Fig:LOREZ:MPC_Comparison}.
The first $10$ time units are used to train with a Schroeder sweep, after which the models are validated on the next $10$ time units using a sinusoidally-based high-frequency forcing, $u(t)=(5\sin(30 t))^3$. 
MPC is then applied for the last $5$ time units.
SINDYc exhibits the best prediction and control performance. 
The neural network exhibits comparable control performance, although the prediction horizon is considerably shorter. % but still sufficient for MPC. 
Surprisingly, DMDc is able to stabilize the fixed point, despite poor predictions based on a linear model.
As the predictive capability of DMDc is poor, we will not present DMDc results in the following, but instead compare SINDYc and the neural network.
As in the previous example, while the neural network exhibits similar control performance, the control execution of SINDYc is $21$ times faster.

\begin{figure}[tb]
\vspace{-.1in}
	\centering
	\includegraphics[width=\textwidth]{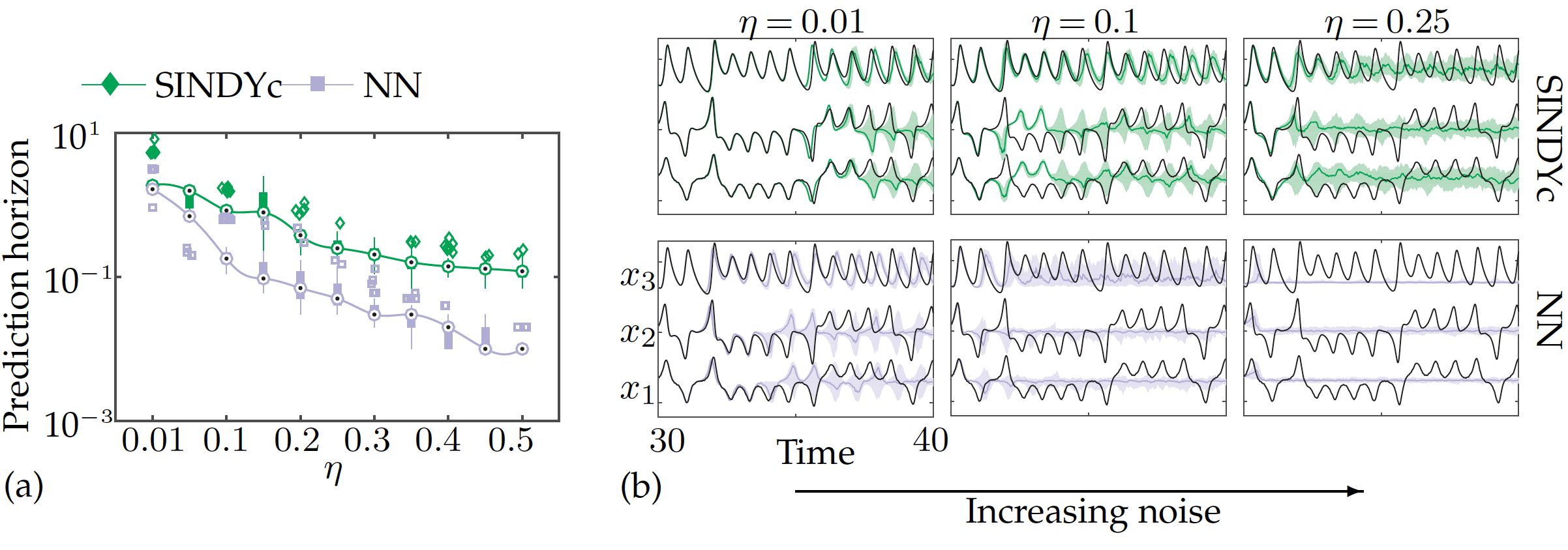}
	\vspace{-.3in}
				\caption{Crossvalidated prediction performance for increasing measurement noise: 
					%			(a) average relative error, 
					(a) prediction horizon in time units, and
					(b) time series with 50 (median as thick colored line) and 25--75 (colored shaded region) percentiles.  
					Statistics are shown for 50 noise realizations each and different noise levels $\eta\in\{ 0.01,0.1,0.25\}$. }
				\label{Fig:LORENZ:NOISE:PredictionPerformance}	
				\vspace{-.3in}
\end{figure}
			
\begin{figure}[tb]
\vspace{-.15in}
	\centering
	\includegraphics[width=\textwidth]{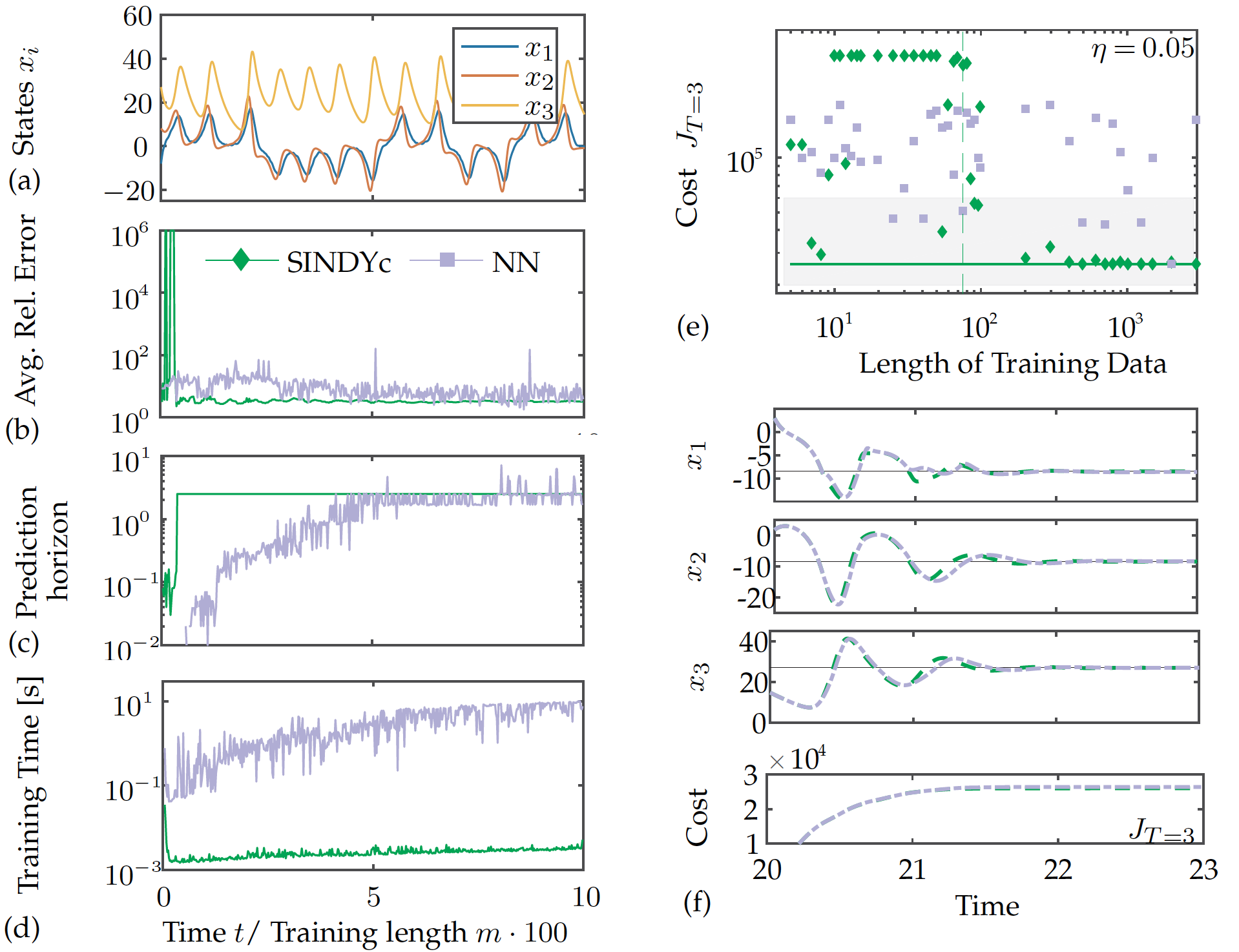}
%	\vspace{-.5in}
\vspace{-.2in}
	\caption{Crossvalidated prediction and control performance for increasing length of training data (without noise): 
		(a) time series of the training data, 
		(b) average relative prediction error, 
		(c) prediction horizon, 
		(d) training time in seconds, 
		(e) terminal cumulative cost over 3 time units, and
		(b) time series of states and cost of the best model for each model type
		($m^{SINDYc} = 38$, $m^{NN} = 40$). 
		Note that from $m = 400$ onwards, SINDY identifies the best performing models.
	}
	\label{Fig:LORENZ:IncreasingTrainingLength}
	\vspace{-.3in}
\end{figure}

Figure~\ref{Fig:LORENZ:NOISE:PredictionPerformance} examines the crossvalidated prediction performance of SINDYc and neural network models based on measurements with increasing noise magnitude $\eta=\sigma/\max(\mathrm{std}(x_i))\in\{0.01,0.1,0.25\}$.  
The performance of both models decreases with increasing noise level, 
although SINDYc generally outperforms the neural network.  
Unlike the Lotka-Volterra model, the average relative error is misleading in this case. 
With increasing noise magnitude the neural network converges to a fixed point, having no predictive power, while SINDYc still exhibits the correct statistics beyond the prediction horizon; however a phase drift leads to a larger average relative error. 
This is shown in Fig.~\ref{Fig:LORENZ:NOISE:PredictionPerformance}(b) with
the median (thick line) and the 25--75 percentile region (shaded area) of the prediction for three different noise levels. 
Thus, a better metric for prediction performance is the prediction horizon itself (see Fig.~\ref{Fig:LORENZ:NOISE:PredictionPerformance}(a)).
The prediction horizon is estimated as the time instant when the error ball is larger than a radius of $\varepsilon = 3$,
i.e.\  a model is considered predictive if $\sqrt{\sum_{i=1}^3 (x_i-\hat{x}_i)^2}<\varepsilon$.
This corresponds to roughly $10\%$ error per state variable, considering that the order of magnitude of each state is approximately $\mathcal{O}(10^1)$; this error radius correlates well with the visible divergence of the true and predicted state in Fig.~\ref{Fig:LORENZ:NOISE:PredictionPerformance}(b).
For low and moderate noise levels, SINDYc robustly predicts the state with high accuracy.
Even for $\eta=0.25$, the 1-period prediction is sufficiently long for a successful stabilization with MPC as we consider a comparably short prediction horizon of $T_p = 0.1$.

The effect of the amount of training data on the prediction and control performance is examined in Figs.~\ref{Fig:LORENZ:IncreasingTrainingLength}, respectively. 
In Fig.~\ref{Fig:LORENZ:IncreasingTrainingLength}(a-d), we show the average relative error evaluated on the prediction over the next $10$ time units, the prediction horizon, and the required training time in seconds for increasing length of noise-free training data.
For a relatively small amount of data, SINDYc rapidly outperforms the neural network model with a prediction horizon of $2.5$ time units and a significantly smaller error.
For a sufficiently large amount of data, SINDYc and the neural network result in comparable predictions. 
However, SINDYc yields highly predictive models that can be rapidly trained in the low and moderate data regimes.
Models trained on weakly noise-corrupted measurements, $\eta=0.05$, are tested in MPC.
For each length of training data, 50 noise realizations are performed and the most predictive model is selected for evaluation in MPC (Fig.~\ref{Fig:LORENZ:IncreasingTrainingLength}(e-f)).
Outside the shaded regions, models are generally not predictive or might even diverge. 
In the noise-corrupted case, it is clear that SINDYc models generally have better control performance than neural network models.
For a sufficiently large amount of training data, neural networks can have comparable performance to SINDYc models, although they show a sensitive dependence on the initial choice of the network weights.
The control results of the neural network are significantly better here than for the Lotka-Volterra model due to the intrinsic system properties.
In chaotic systems, a long enough trajectory will come arbitrarily close to every point on the attractor; thus, measurements of the Lorenz system are in some sense richer than those of the Lotka-Volterra system.
A surprising result is that a nearly optimal SINDYc model can be trained on just $8$ noisy measurements (compare Fig.~\ref{Fig:LORENZ:IncreasingTrainingLength}(e-f)).

%% file: Sec43.tex
\begin{figure}[tb]
	\centering
	\includegraphics[width=\textwidth]{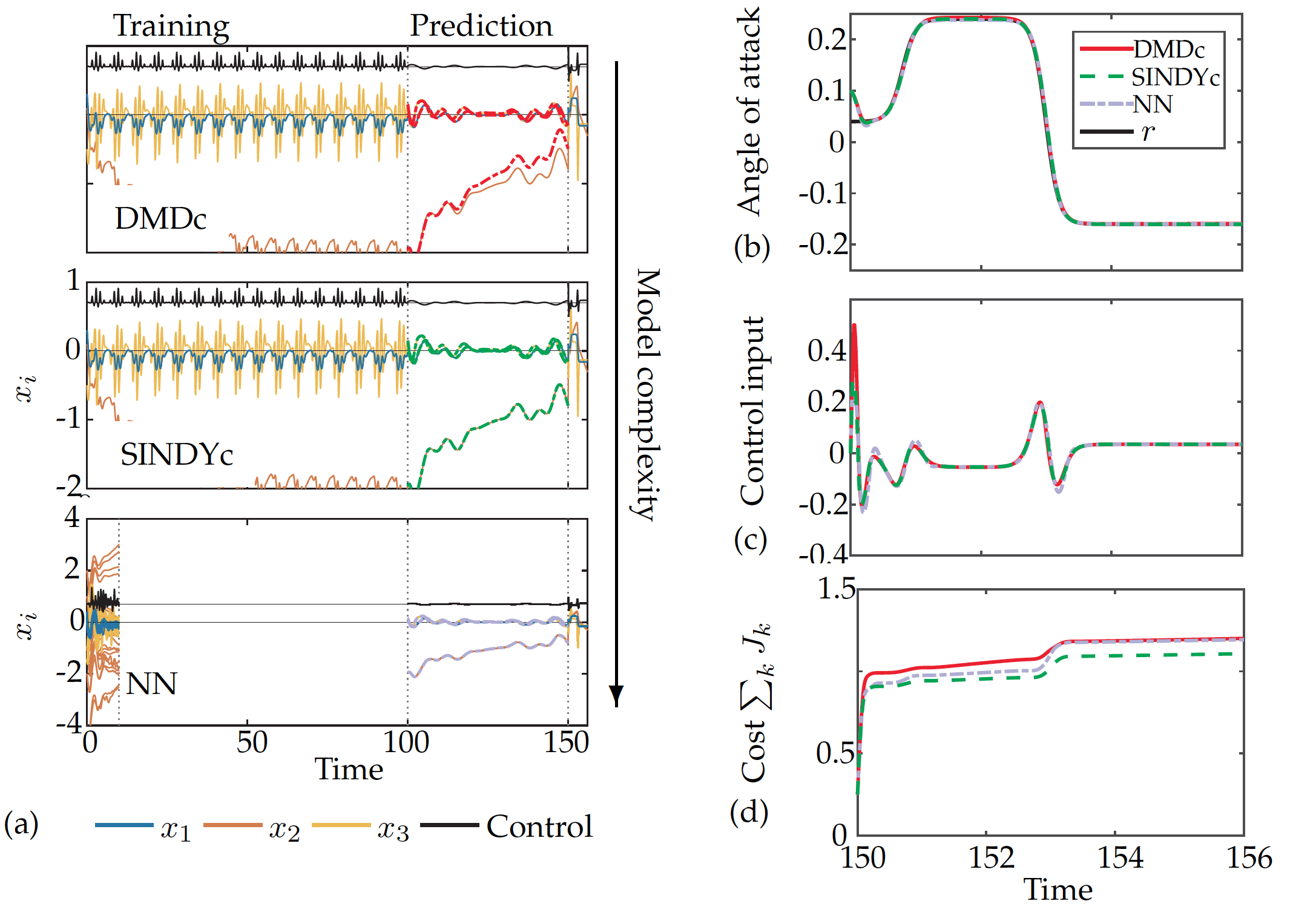}
	\vspace{-.3in}
	\caption{Prediction and control performance for F8 aircraft:
		(a) time series of states and input during training, validation, and control, 
		(b) angle of attack $y=x_1$ with reference $r$,
		(c) control input, and 
		(d) cumulative cost during control.
		%The MPC optimization requires about 2min for DMDc, 6min for SINDYc and 70min for the neural network.
	}
	\label{Fig:F8:MPC_Comparison}
	\vspace{-.3in}
\end{figure}

\section{Tracking for the F-8 crusader}\label{Sec:F8-Aircraft}
In this section we consider an automatic flight control system of the F-8 aircraft at an altitude of 30 000 ft (9000 m) and Mach=0.85~\cite{garrard1977automatica,cimen2004scl,yan2012ieee}. The control objective is to track a specific trajectory of the angle of attack. The aircraft dynamics~\cite{garrard1977automatica} are given by
\begin{subequations}\label{Eqn:F8}
	\begin{align}
		\dot{x}_1 &= -0.877x_1 +x_3 -0.088x_1 x_3 +0.47x_1^2 -0.019 x_2^2 -x_1^2 x_3 + 3.846 x_1^3 - 0.215 u\\ 
					  &\quad + 0.28 x_1^2 u + 0.47x_1 u^2 + 0.63 u^3\notag \\ 
		\dot{x}_2 &= x_3 \\
		\dot{x}_3 &= -4.208 x_1 - 0.396 x_3 - 0.47 x_1^2 - 3.564 x_1^3 - 20.967 u + 6.265 x_1^2 u   + 46 x_1 u^2 + 61.1 u^3
	\end{align}
\end{subequations}
where $x_1$ is the angle of attack (rad), $x_2$ is the pitch angle (rad), $x_3$ is the pitch rate (rad s$^{-1}$), and $u$ is the control input representing the tail deflection angle (rad).
The system is nonaffine in the states and the control input rendering it strongly nonlinear. 
The commanded angle of attack to be tracked~\cite{yan2012ieee} is given by:
\begin{equation}
r(t) = 0.4 \left( - \frac{0.5}{1+e^{\hat{t}-0.8}} + \frac{1}{1+e^{\hat{t}-3}}-0.4\right)
\end{equation}
with $\hat{t} = t/0.1$.
We assume that the output, over which the performance is optimized, is $y = x_1$.
The timestep of the system is $\Delta t = 0.001$ and the timestep of the model is $\Delta t^{M} = 0.01$. The control input is determined using SINDY-MPC every $10$ system timesteps over which the applied control is then kept constant.
The weight matrices are 
$Q=25$,
%$Q = \left(\begin{smallmatrix} 25 & 0 & 0\\ 0 & 0 & 0\\ 0 & 0 & 0\end{smallmatrix}\right)$, 
$R_u =  R_{\Delta u} = 0.05$, the actuation input rate is limited to $\Delta u\in[-0.3,0.5]$, and the constraint for the angle of attack is $y\in[-0.2,0.4]$.
The control and prediction horizon is $m_p = m_c = 13$ and the sparsity-promoting parameter in SINDYc is $\boldsymbol{\lambda} = (10^{-4},10^{-2},10^{-2})$, where $\lambda_i$ is used to identify the terms for $x_i$. 
The neural network has two hidden layers each with 15 neurons. 
Access to full-state information is assumed for these models.

Results assessing prediction and control performance of SINDYc compared with DMDc and a neural network model are displayed in Fig.~\ref{Fig:F8:MPC_Comparison}. 
Similar to the Lotka-Volterra system, the neural network requires more and {\it richer} training data, i.e. a better exploration of the system behavior, to perform sufficiently well. 
Thus, 250 short trajectories each consisting of $1000$ snapshots ($25\cdot 10^4$ instances in total) with varying input signals are used to train the neural network; a subset of 20 trajectories is displayed in Fig.~\ref{Fig:F8:MPC_Comparison}(a,bottom).
In contrast, SINDYc and DMDc perform similarly well if trained on much less data ($10^4$ instances of a single trajectory). 
Moreover, SINDYc learns from few measurements the true relationship between the variables, even though only limited system behavior has been observed, resulting in increased performance.

%% file: Sec44.tex
%\vspace{-1in}
\section{Optimal therapy for pathogenic attacks}\label{Sec:HIV}
Optimizing drug therapy is critical for inhibiting diseases such as cancer and viral infections.
Here, we consider treatment of infections with the human immunodeficiency virus (HIV), 
a pathogen that infects T-helper CD4+ cells of the immune system and can cause Acquired Immune Deficiency Syndrome (AIDS).
Identifying the underlying infection mechanism, the response of the immune system, and the interactions with drugs targeting different components in this system is critical for developing and optimizing therapeutic strategies.
Various models have been proposed to study the interaction between HIV and CD4+ cells; we refer to a recent review~\cite{perelson1999siam}. 

Optimal treatment aims to decrease virus mutations, complications from administered drugs, medical costs, and to strengthen the immune system.
We consider a system~\cite{zurakowski2006jtb} that incorporates infections with HIV, the cytotoxic lymphocyte (CTL) response of the immune system, and therapeutic interventions via a highly active anti-retroviral therapy (HAART), i.e.\ a combination of drugs that affect the replication rate of HIV and support the immune system. 
%This is based on a more general system that is simplified under certain conditions~\cite{wodarz2001jtb}:
This is based on a more general and complex system, that can be simplified under certain conditions~\cite{wodarz2001jtb}:
% blocking new infections
%
\begin{subequations}\label{Eqn:HIV}
	\begin{align}
	\dot{x}_1 &= \lambda - d x_1 - \beta(1-\eta u)x_1 x_2\\ 
	\dot{x}_2 &= \beta(1-\eta u)x_1 x_2 - a x_2 -p_1x_4x_2 - p_2x_5x_2\\
	\dot{x}_3 &= c_2x_1x_2x_3 - c_2qx_2x_3-b_2x_3\\
	\dot{x}_4 &= c_1x_2x_4 - b_1x_4\\
	\dot{x}_5 &= c_2qx_2x_3 - hx_5
	\end{align}
\end{subequations}
with parameters $\lambda=1$, $d=0.1$, $\beta=1$, $a=0.2$, $p_1 = 1$, $p_2 = 1$, $c_1 = 0.03$, $c_2 = 0.06$, $b_1 = 0.1$, $b_2 = 0.01$, $q=0.5$, $h=0.1$ and $\eta = 0.9799$ (units typically in mm$^{-3}$/day).
Here, the states describe concentrations of healthy CD4+ T-cells, $x_1$, HIV-infected CD4+ T-cells, $x_2$ ,
CTL precursors (memory CTL), $x_3$,
helper-independent CTL, $x_4$, and helper-dependent CTL, $x_5$.
For a detailed discussion of the system~\eqref{Eqn:HIV} we refer to~\cite{zurakowski2006jtb,wodarz2001jtb}.
The parameter $\eta$ represents the effectiveness of the HAART therapy applied via $u$. %, where $0\leq u\leq1$.
%Healthy cells $x_1$ are produced at rate $\lambda$ and die at rate $dx_1$, 
%$k_1x_1x_j$ for $j=4,5$ models infection process:
%cells $x_2$ are produced by unaffected T-cells $x_1$ and virus $x_j$, $4,5$.
%drug-sensitive wild-type and mutant viruses infect target cells by rates $k_1$ and $k_2$, respectively, and rates at which infected cells produce the virus, $\alpha_1$ and $\alpha_2$, respectively.
%Infected cells and viruses die at rates $\delta x_j$ with $j=2,3$ and $cV_j$ with $j=4,5$.
%probability of mutation from drug-sensitive wide-type virus to mutant virus is given by $\mu$. 
For the considered parameters and in the absence of control ($u\equiv 0$), the system exhibits two stable fixed points: a progressive infection leading to AIDS, $\bx^A$, and
the recovery from a successful immune response, $\bx^B$.
The later steady-state is given by
\begin{subequations}\label{Eqn:HIV_SSB}
	\begin{align}
		x^B_1 &= \frac{\lambda}{d+\beta x^B_2},\quad
		x^B_2 = \frac{c_2(\lambda-dq)-b_2\beta - \sqrt{[c_2(\lambda-dq)-b_2\beta]^2-4\beta c_2qdb_2}}{2\beta c_2q},\\
		x^B_3 &= \frac{hx^B_5}{c_2 qx^B_2},\quad
		x^B_4 = 0,\quad
		x^B_5 = \frac{x^B_2 c_2(\beta q-a)+b_2\beta}{c_2p_2x^B_2},
	\end{align}
\end{subequations}
and exists if $[c_2(\lambda-dq)-b_2\beta]^2-4\beta c_2qdb_2 \geq0$.
The region of attraction (ROA) to this fixed point is limited and only established if the infectivity of the virus is small such that 
$\beta< \frac{c_1[c_2b_2(\lambda-qd)-b_2c_1d]}{b_1(c_2b_1q+b_2c_1)}$, which can be achieved by applying a HAART therapy ($u=1$) with high efficacy ($\eta\approx 1$). 
%This state moves as a function of $\beta_{eff} = \beta(1-\eta u)$ when $u=1$ and its ROA lies outside of the ROA in the absence of drug treatment ($u=0$), so the control policy must first drive the system close to the desired fixed point $\bx^B$. 
This state moves as a function of $\beta_{eff} = \beta(1-\eta u)$ when $u>0$ and its ROA changes and does not necessarily overlap with the ROA in the absence of drug treatment ($u=0$), i.e.\ dependent on the initial condition and the applied control the system will converge to a different steady-state. 
A non-trivial control strategy is required that switches between treatment and no treatment to establish a successful immune response, and hence to approach $\bx_0^B$. 
In contrast, when treatment is applied continuously for a sufficiently long amount of time such the fixed points are approached and then terminated, the system will converge to a progressive infection, $\bx^A$, even if a successful immune response had been established.

\begin{figure}[tb]
	\centering
	\includegraphics[width=\textwidth]{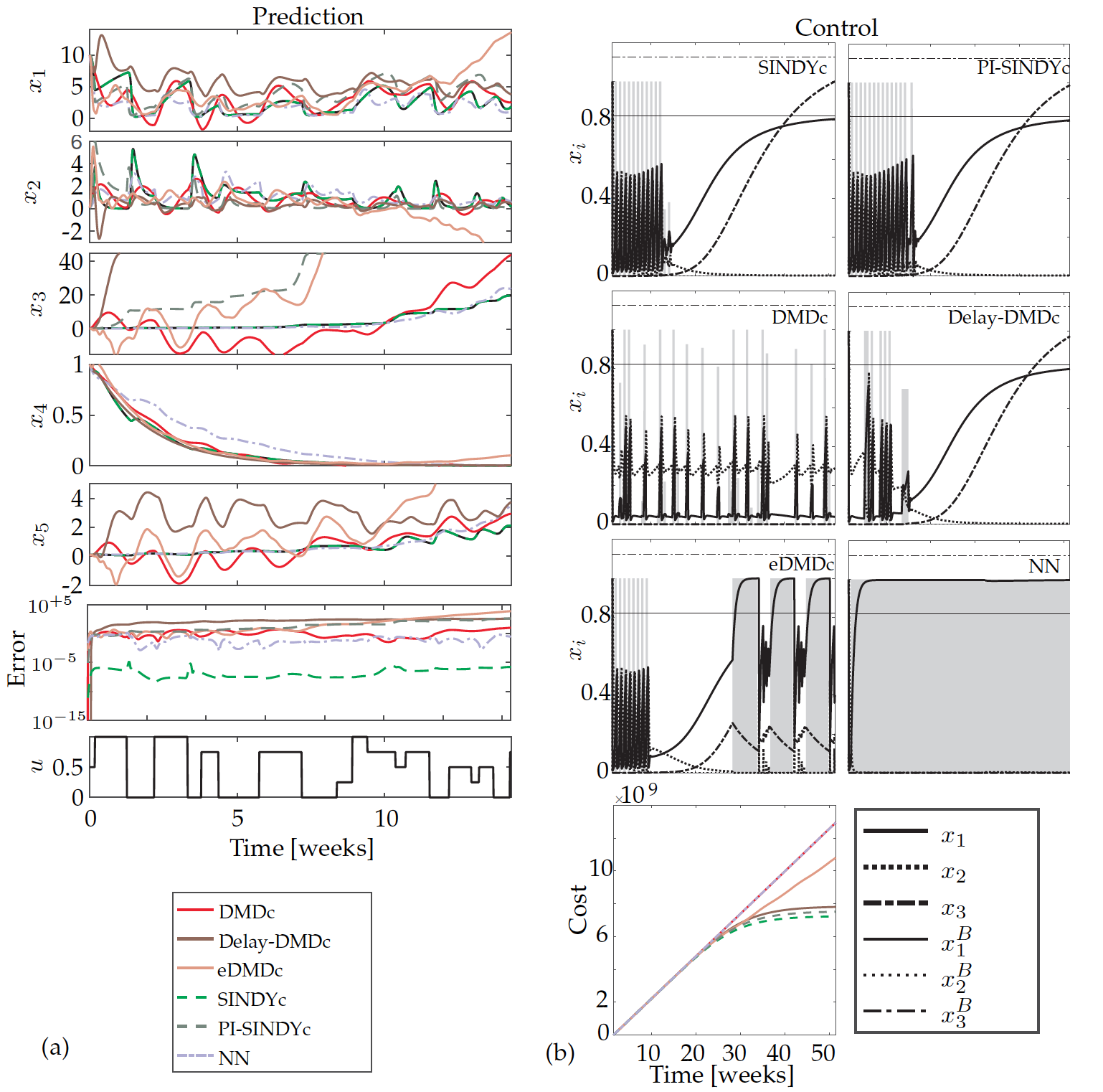}
	\vspace{-.3in}
	\caption{Prediction and control performance for various models of the HIV model: 
		(a) Predicted and true (solid black) time series of the forced dynamics, and
		(b) control results (normalized for better legibility) with optimized actuation input (shaded background) and reference values. 
		Only SINDYc and, with slightly less performance, PI-SINDYc and Delay-DMDc are capable of driving the system to the desired steady state.}
	\label{Fig:HIV:MPC_Comparison}
	\vspace{-.2in}
\end{figure}

The cost functional to be optimized is given by
\begin{equation}
J = \int_{0}^{T}\, 
(x_1(t)-\hat{x}_1) + (x_3(t)-\hat{x}_3) + \vert u(t)\vert \,dt
\end{equation}
where $\hat{x}_1 = x^B_1$ and $\hat{x}_3 = x^B_3$~\cite{zurakowski2006jtb}
taking into account the healthy cells, the immune system, and the cost of treatment.
The control input $u$ is bounded by $0\leq u\leq 1$ with efficacy of $\eta = 0.9799$.
An additional constraint is added to the control that renders all cell concentrations nonnegative, i.e.\ $x_i\geq 0\;\forall i$.
The time step is $\Delta t^M = 2\,\mathrm{hrs}$ for the model and  is $\Delta t = 1/24\,\mathrm{day} =  1\,\mathrm{hr}$ for the simulated system. 
The control performance is evaluated over $50$ weeks. The prediction and control horizon for the MPC optimization are both $m_p=m_c=24$, i.e.\ over 2 days (from $m_p\Delta t^M$). 
We assume a more realistic situation, where the state is measured once a week, and the treatment is then kept constant over the following week.
The training data consists of samples collected over $200$ days ($\approx 30$ weeks) with a discrete control input, as was applied for the validation data in Fig.~\ref{Fig:HIV:MPC_Comparison} (bottom). 
In contrast, the training data for the neural network consists of ensemble data of $32$ different trajectories. 
In both cases, the control is a random sequence of values that are kept constant over random durations of time $\Delta T\in[5\,\mathrm{hrs},10\,\mathrm{days}]$. 
 
 \begin{table*}[tb]
 	\begin{center}
 		\begin{tabular*}{\textwidth}{V{0.15\textwidth} V{0.3\textwidth} V{0.35\textwidth} V{0.2\textwidth}} %3.5
 			\toprule
 			Model & State & Settings & \# Unknowns\\
 			\midrule
 			SINDYc & $\by = \bx$ & polynomial basis (order $r=3$), $\boldsymbol{\lambda} = (10,3.1,3,0.1,0.5)$& $84\cdot 5 = 420$ \\ 
 			PI-SINDYc & $\by = (x_1,x_2,x_3)$ & polynomial basis (order $r=3$), $\boldsymbol{\lambda} = (10,30,3)$ & $35\cdot 3 = 105$\\
 			DMDc & $\by =\bx-\bx^B$ & Deviation from reference state & $5^2=25$\\
 			Delay-DMDc & $\by =(\bx(t)-\bx^B,\bx(t-\tau)-\bx^B,\ldots,\bx(t-9\tau)-\bx^B)$ &  Deviation from reference state, 10 time delay coordinates of the full-state and of the control input& $(10\cdot5)^2=2500$\\
 			eDMDc & $\by =\boldsymbol{\Theta}(\bx,r)$ & order $r=3$ of polynomial basis (without constant term) & $83^2 = 6889$\\
 			NN & $\by =\bx$ & 1 hidden layer with 5 neurons and linear activation functions, data is $log$-transformed and mapped to $[-1,\; 1]$ to compensate for skewness and different range & 43\\
 			\bottomrule	
 		\end{tabular*}
		\vspace{-.1in}
 		\caption{Model parameters for the HIV system.}
 		\label{Tab:HIV:ModelParameters}
 	\end{center}
 	\vspace{-.45in}
 \end{table*}	
We consider a SINDYc model with (1) full-state information (SINDYc) and (2) partial information based on a subset of the variables (PI-SINDYc).
The latter case demonstrates the situation when only a few states can be measured, which is generally more realistic. 
For the identification of the SINDYc models it is important to normalize first the features in the library, as the coefficients of the active terms spread over several orders of magnitude.
In both cases a polynomial order of three is used for the library.
The results are compared with various linear models: (1) DMDc on the full state (DMDc), (2) DMDc on delay coordinates of the full state (Delay-DMDc), and (3) DMDc on a set of nonlinear observables (extended DMD with control, eDMDc).
In addition, a neural network model on the full state (NN) is trained. An overview of these models and their parameters is provided in Tab.~\ref{Tab:HIV:ModelParameters}.
The identified parameters of the SINDYc and PI-SINDYc models are displayed below~\eqref{Eqn:HIV:SINDYc} are:
\begin{eqnarray}\label{Eqn:HIV:SINDYc}
\scriptsize
%\boldsymbol{\Xi}^{SINDYc} &=& \scriptsize
\underbrace{\left[\begin{array}{ccccc}
	\texttt{0.9995} & \texttt{0} & \texttt{0} & \texttt{0} & \texttt{0}\\  %1
	\texttt{-0.0999} & \texttt{0} & \texttt{0} & \texttt{0} & \texttt{0}\\  %x1
	\texttt{0} & \texttt{-0.1991} & \texttt{0} & \texttt{0} & \texttt{0}\\  %x2
	\texttt{0} & \texttt{0} & \texttt{-0.0100} & \texttt{0} & \texttt{0}\\  %x3
	\texttt{0} & \texttt{0} & \texttt{0}  & \texttt{-0.1000} & \texttt{0}\\ %x4
	\texttt{0} & \texttt{0} & \texttt{0}   & \texttt{0} & \texttt{-0.1000}\\%x5
	\texttt{-0.9990} & \texttt{0.9981} & \texttt{0}  & \texttt{0} & \texttt{0}\\ %x1x2
	\texttt{0} & \texttt{0} & \texttt{-0.0299} & \texttt{0} & \texttt{0.0300}\\  %x2x3
	\texttt{0} & \texttt{-0.9982} & \texttt{0}  & \texttt{0.0300} & \texttt{0}\\ %x2x4
	\texttt{0} & \texttt{-0.9990} & \texttt{0} & \texttt{0} & \texttt{0}\\  %x2x5
	\texttt{0} & \texttt{0} & \texttt{0.0600} & \texttt{0} & \texttt{0}\\  %x1x2x3
	\texttt{0} & \texttt{0} & \texttt{0}  & \texttt{0} & \texttt{0}\\ %x1x2x4
	\texttt{0} & \texttt{0} & \texttt{0} & \texttt{0} & \texttt{0}\\  %x1x2x5
	\texttt{0.9763} & \texttt{-0.9757} & \texttt{0} & \texttt{0} & \texttt{0}  %x1x2u
	\end{array}\right]}_{\boldsymbol{\Xi}^{SINDYc} }\quad
%\boldsymbol{\Xi}^{\scriptsize PI-SINDYc} = \scriptsize
\underbrace{\left[\begin{array}{ccc}
	\texttt{0.9995} & \texttt{0} & \texttt{0}\\  %1
	\texttt{-0.0999} & \texttt{0} & \texttt{0}\\  %x1
	\texttt{0} & \texttt{-0.8985} & \texttt{0}\\  %x2
	\texttt{0} & \texttt{0} & \texttt{-0.0100}\\   %x3
	\texttt{0} & \texttt{0} & \texttt{0}\\   %x4
	\texttt{0} & \texttt{0} & \texttt{0} \\  %x5
	\texttt{-0.9990} & \texttt{0.9424} & \texttt{0} \\  %x1x2
	\texttt{0} & \texttt{-0.0573} & \texttt{-0.0299}\\   %x2x3
	\texttt{0} & \texttt{0} & \texttt{0} \\  %x2x4
	\texttt{0} & \texttt{0} & \texttt{0} \\  %x2x5
	\texttt{0} & \texttt{0.0069} & \texttt{0.0600}  \\ %x1x2x3
	\texttt{0} & \texttt{0} & \texttt{0}\\   %x1x2x4
	\texttt{0} & \texttt{0} & \texttt{0} \\  %x1x2x5
	\texttt{0.9763} & \texttt{-0.7507} & \texttt{0}   %x1x2u
	\end{array}\right]}_{\boldsymbol{\Xi}^{PI-SINDYc}} 
\begin{array}{c}
1\\
x_1\\
x_2\\
x_3\\
x_4\\
x_5\\
x_1x_2\\
x_2x_3\\
x_2x_4\\
x_2x_5\\
x_1x_2x_3\\
x_1x_2x_4\\
x_1x_2x_5\\
x_1x_2 u
\end{array}
\end{eqnarray}
Only the non-zero parameters are shown, and their error is $\mathcal{O}(10^{-3})-\mathcal{O}(10^{-6})$ for SINDYc.
The error in the parameters decreases with increased time resolution.
Here, a coarse time step is chosen to reduce the computational cost of MPC for the chosen prediction horizon.
In PI-SINDYc, the parameters for $x_1$ and $x_3$ are estimated well as these only depend on $x_1$, $x_2$ and $x_3$. In contrast, $x_2$ has a larger error in the estimated parameters and consists of erroneous parameters to compensate for the missing information. Different selections of variables have been tested, which generally resulted in poor models, except for the selected combination.
The resulting models are generally not sparse, except where a direct relationship exists between variables. This suggests that SINDY indicates direct causal relationships, which can be measured in terms of the sparsity.  

Prediction accuracy based on data differing from the training set, but with a similar type of actuation signal, and control results are displayed in Fig.~\ref{Fig:HIV:MPC_Comparison}.
Both start from an early infection given by $\bx_0 = (\lambda/d, 0.1, 0.1, 0.1, 0.1)^T$.
While a SINDYc model can be identified with near-perfect prediction accuracy, all other models display an error several orders of magnitude larger (see Fig.~\ref{Fig:HIV:MPC_Comparison}(a)). 
In particular, linear DMDc-based models diverge significantly from the true trajectory for some variables, while capturing the right trend in other variables.
The neural network and the PI-SINDYc model based on partial state information generally stay closer to, and even temporarily match, the true trajectory.
Interestingly, while MPC using PI-SINDYc successfully drives the system to the desired steady-state behavior, with a slightly larger cost than SINDYc, 
the neural network controller is unable to establish the successful immune response by applying constant treatment (see~Fig.~\ref{Fig:HIV:MPC_Comparison}(b)). 		
Note that the actuation depends strongly on the prediction and control horizon chosen for the optimization; further analysis has shown that a smaller horizon for the NN controller yields a time-varying, however still unsuccessful, treatment.
We varied the number of hidden layers (up to 3), the number of neurons (up to 100), the type of activation function, the number of delays (up t0 100) in the state and input variables and the amount of training data ($\approx 600$ different initial conditions). However, these did not significantly change the performance of the model. The type of data (not just the amount) is particularly critical for training a neural network. Designing experiments, i.e.\ a good forcing signal that explores the system behavior and yields \textit{dynamically rich} training data, is a challenge of its own.   
The linear DMDc and eDMDc models fail too. While the eDMDc model starts with the correct frequency, detrimental treatment is administered thereafter close to the desired state, which gives rise to new growth of infected cells, $x_2$. Interestingly, augmenting the state vector with delay coordinates results in a successful treatment (with performance close to the SINDYc models), in contrast to the strategy to augment the state with nonlinear measurements of the state as in eDMDc.

All models but the neural network, which has been trained on a significantly larger amount of data, have been trained on the same amount of data, a single trajectory starting from an initial condition which is relatively far from the desired behavior.
Thus, these models are required to generalize well, i.e.\ perform well far from the region in which they have been initially trained.
Using more data would certainly help to improve the prediction accuracy of some of these models, in particular, if these require a large number of parameters to be estimated. 
However, this would pose additional challenges in real-time applications with abrupt system changes, as this requires robust model formation and adaptation from few measurements.

%% file: Sec5.tex
\input{Sec5Table}

\section{Discussion and Conclusions}\label{Sec:Discussion}
In conclusion, we have demonstrated the effective integration of data-driven sparse model discovery for model predictive control in the low-data limit.  
The sparse identification of nonlinear dynamics (SINDY) algorithm has been extended to discover nonlinear models with actuation and control, resulting in interpretable and parsimonious models.  
Moreover, because SINDY only identifies the few active terms in the dynamics, it requires less data than many other leading machine learning techniques, such as neural networks, and prevents overfitting.  
When integrated with model predictive control, SINDY provides computationally tractable and accurate models that can be trained on very little data.  
The resulting SINDY-MPC framework is capable of controlling strongly nonlinear systems, purely from measurement data, and the model identification is fast enough to discover models in real-time, even in response to abrupt changes to the model.   
The SINDY-MPC approach is compared with MPC based on data-driven linear models and neural network models on four nonlinear dynamical systems of different complexities and challenges: the weakly nonlinear Lotka-Volterra system, the chaotic Lorenz system, the nonaffine F8 crusador model, and the HIV/immune response system, which variables are of different order of magnitudes and where only partial state information is available. 

The relative strengths and weaknesses of each method are summarized in Tab.~\ref{Tab:Summary}.  
By nearly every metric, linear DMDc models and nonlinear SINDYc models outperform neural network models (NN).  
In fact, DMDc may be seen as the limit of SINDYc when the library of candidate terms is restricted to linear terms.  
SINDY-MPC provides the highest performance control and requires significantly less training data and execution time compared with NN.  
However, for very low amounts of training data, DMDc provides a useful model until the SINDYc algorithm has enough data to characterize the dynamics.  
Thus, we advocate the SINDY-MPC framework for effective and efficient nonlinear control, with DMDc as a stopgap after abrupt changes until a new SINDYc model can be identified.  
\RefOne{
Note that a crucial step in SINDY is the choice of library functions, which is often informed by expert knowledge about what category of nonlinearities to include. 
A poor choice of the library will generally yield a non-sparse model. 
Without any prior knowledge about the system type, a sweep through different classes of candidate functions is required. 
However, once a model is learned from a sufficiently rich library, the model is often able to generalize beyond the training data.  
%If we don't know anything about the system, a sweep through classes of candidate functions would be required. 
%However, once having learned a model from a sufficiently rich library, the model is able to generalize beyond the training data.
If the model structure is not fixed, but varies heterogeneously in state space, neural networks may provide a more flexible and generalizable architecture to represent the dynamics.  % provide generalize better on regions the model has not been trained on.
A heterogeneous model structure can potentially be incorporated into SINDy by additionally learning a library of models~\cite{Brunton2014siads,sargsyan2015pre}.
}

This work motivates a number of future extensions and investigations.  
Although the preliminary application of SINDYc for MPC is encouraging, this study does not leverage many of the powerful new techniques in sparse model identification.  
Figure~\ref{Fig:SINDYc:Schematic} provides a schematic of the modularity and demonstrated extensions that are possible within the SINDy framework.  
In realistic applications, the system may be extremely high-dimensional, and the SINDy library does not scale well with the size of the data.  
Fortunately, many high-dimensional systems evolve on a low-dimensional attractor, and it is often possible to identify a model on this attractor, for example by identifying a SINDy model on low-dimensional coordinates obtained through a singular value decomposition~\cite{Brunton2016pnas} or manifold learning~\cite{Loiseau2018jfm}.  
In other applications, full-state measurements are unavailable, and the system must be characterized by limited measurements.  
It has recently been shown that delay coordinates provide a useful embedding to identify simple models of chaotic systems~\cite{Brunton2017natcomm}, building on the celebrated Takens embedding theorem~\cite{Takens1981lnm}.  
Delay coordinates also define intrinsic coordinates for the Koopman operator~\cite{Brunton2017natcomm}, which provides a simple linear embedding of nonlinear systems~\cite{Mezic2005nd,Budivsic2012chaos}.  
Koopman models have recently been used for MPC~\cite{korda2016_a,Peitz2017arxiv} and have been identified using SINDy regression~\cite{Kaiser2017arxiv} and subsequently used for optimal control~\cite{Kaiser2017arxiv}.  
Recently, SINDY has been extended to modify an existing model based on new incoming measurements to enable rapid model recovery from abrupt changes to the system~\cite{Quade2018arxiv}.
\RefTwo{
Learning quickly from limited measurements is an important task, which may be viewed in terms of design of experiments; specifically, optimizing the actuation input to collect the most informative measurements to learn a more predictive model faster. 
This would require the formulation of a different cost function, which measures the predictive power of the model, to determine future actuation inputs.
Rapid learning is also related to the question of quantity versus quality of data and identifiability~\cite{gevers2013cdc,alkhoury2017automatica}; more data is usually better, although it is possible to work with less data if it is representative of the system.  
}
Further, similar methods could be used to optimize sensors and exploit partial measurements within the SINDY-MPC framework.  
All of these innovations suggest a shift from the perspective of \emph{big data} to the control-oriented perspective of \emph{smart data}.  

Figure~\ref{Fig:SINDYc:Schematic} also demonstrates innovations to the SINDy regression to include physical constraints, known model structure, and model selection, which may all benefit the goal of real-time identification and control.  
Known symmetries, conservation laws, and constraints may be readily included in both the SINDYc and DMDc modeling frameworks~\cite{Loiseau2016jfm}, as they are both based on least-squares regression, possibly with sequential thresholding.  
It is thus possible to use a constrained least-squares algorithm, for example to enforce energy conserving constraints in a fluid system, which manifest as anti-symmetric quadratic terms~\cite{Loiseau2016jfm}.  
Enforcing constraints has the potential to further reduce the amount of data required to identify models, as there are less free parameters to estimate, and the resulting systems have been shown to have improved stability in some cases.  
It is also possible to extend the SINDy algorithm to identify models in libraries that encode richer dynamics, such as rational function nonlinearities~\cite{Mangan2016ieee}.  
Finally, incorporating information criteria provides an objective metric for model selection among various candidate SINDy models with a range of complexity.  

The SINDY-MPC framework has significant potential for the real-time control of strongly nonlinear systems.  
Moreover, the rapid training and execution times indicate that SINDy models may be useful for rapid model identification in response to abrupt model changes, and this warrants further investigation. 
The ability to identify accurate and efficient models with small amounts of training data may be a key enabler of recovery in time-critical scenarios, such as model changes that lead to instability.  
In addition, for broad applicability and adoption, the SINDy modeling framework must be further investigated to characterize the effect of noise, derive error estimates, and provide conditions and guarantees of convergence.  
These future theoretical and analytical extensions are necessary to certify the model-based control performance.

%% file: Sec5Table.tex
\begin{table*}[tb]
\small
	\begin{center}
	\begin{tabular*}{\textwidth}{V{2cm} V{3.5cm} V{3.5cm} V{3.5cm}}
		\toprule
		Property & \textbf{DMDc} & \textbf{SINDYc} & \textbf{NN}\\
		\midrule
		Training with  &\textbf{\textcolor{ForestGreen}{\underline{strong}}}& \textcolor{ForestGreen}{\bf strong} & \textcolor{red}{\bf weak}\\[-0.5ex]
		limited data	& \small Very few samples are sufficient. & \small Well suited for low and medium amount of data.  & \small Requires long time series to learn predictive models.\\[1ex]
		High-  &\textbf{\textcolor{ForestGreen}{\underline{strong}}}& \textcolor{orange}{\bf fair} & \textcolor{ForestGreen}{\bf strong}\\[-0.5ex]							   
		dimensionality & \small Can handle high-dim. data in combination with SVD. & \small Limited by the library size. & \\[1ex]
		Nonlinearities & \textbf{\textcolor{red}{weak}}/\textcolor{orange}{fair} & \textbf{\textcolor{ForestGreen}{strong}}  & \textbf{\textcolor{ForestGreen}{strong}} \\[-0.5ex]
		& \small Linear and weakly nonlinear, however with performance loss. &  \small  Suitable for strongly nonlinear systems. & \small Suitable for strongly nonlinear systems.\\[1ex]						
		Prediction performance & \textbf{\textcolor{orange}{fair}} & \textbf{\textcolor{ForestGreen}{\underline{strong}}} & \textbf{\textcolor{ForestGreen}{strong}}\\[-0.5ex]
		& & & \\[-1ex] % [1ex]
		Control performance  & \textbf{\textcolor{orange}{fair}} & \textbf{\textcolor{ForestGreen}{\underline{strong}}} & \textbf{\textcolor{ForestGreen}{strong}}\\[-0.5ex]
		& & & \\[-1ex] % [1ex]
		Noise  & \textbf{\textcolor{red}{weak}} & \textbf{\textcolor{ForestGreen}{\underline{strong}}} & \textbf{\textcolor{orange}{fair}}\\[-0.5ex]
		robustness & \small High sensitivity w.r.t. noise.& \small Intrinsic robustness due to sparse regression. & \small Can handle low noise levels.\\[1ex]								
		Parameter  & \textbf{\textcolor{ForestGreen}{\underline{strong}}} & \textbf{\textcolor{ForestGreen}{strong}} & \textbf{\textcolor{red}{weak}}\\[-0.5ex]
		robustness & & & \small High sensitivity w.r.t. initial weights of the network.\\[1ex]		
		Training time & \textbf{\textcolor{ForestGreen}{\underline{strong}}}  & \textbf{\textcolor{ForestGreen}{strong}}  & \textbf{\textcolor{red}{weak}}\\[-0.5ex]
& & & \\[-1ex] % [1ex]					
		Execution  & \textbf{\textcolor{ForestGreen}{\underline{strong}}}  & \textbf{\textcolor{ForestGreen}{strong}}  & \textbf{\textcolor{red}{weak}}\\[-0.5ex]
		time & \small Fast optimization routines exist for linear systems. & & \\[1ex]		
		\bottomrule
	\end{tabular*}
	\vspace{-.1in}
	\caption{Capabilities and challenges of DMDc, SINDYc and NN  models. The model with the \textcolor{ForestGreen}{\underline{strongest}} performance is underlined.}\label{Tab:Summary}
		\end{center}
			\vspace{-0.5in}
\end{table*}